\documentclass[a4paper,11pt]{article}

\usepackage{amsmath,amsfonts,amssymb,amsthm}
\usepackage{mathtools}
\usepackage{stmaryrd}
\usepackage{latexsym}
\usepackage{cite}
\usepackage[dvips]{graphicx}
\usepackage[small,bf]{caption}
\usepackage{multirow}
\usepackage{float}
\usepackage{subfig}

\usepackage{listings}
\usepackage{xcolor}
\lstdefinelanguage{pseudo}{
frame=BT,
mathescape,
morekeywords={1.,2.,3.,4.,5.,6.,7.,8.,9.,10.,11.,12.,(a),(b),(c),(d)},
basicstyle=\normalsize \ttfamily \color{black},
showspaces=false, 
showstringspaces=false,        
showtabs=false,
keywordstyle=\bfseries,
}

\newcommand\mytop[2]{\genfrac{}{}{0pt}{}{#1}{#2}}
\newcommand{\enorm}[1]{\,|\!|\!| {{#1}} |\!|\!|_{\Omega}}

\def \Id {\operatorname{Id}}
\def \Hdiv {H(\operatorname{div})}
\def \etal {\emph{et al.}}
\numberwithin{equation}{section}
\theoremstyle{plain}
\newtheorem{theorem}{Theorem}[section]
\theoremstyle{remark}
\newtheorem{rmrk}[theorem]{Remark}
\usepackage[hang,flushmargin]{footmisc}
\providecommand{\keywords}[1]{{\small{\textbf{Keywords:}} #1}} 

\title{Certified Descent 
Algorithm for shape optimization driven by fully-computable a posteriori error 
estimators
}
\author{
M. Giacomini \footnotemark[1]\textsuperscript{ \ ,}\footnotemark[2] , O. Pantz  
\footnotemark[3] \ and K. Trabelsi \footnotemark[2]
}
\date{}

\begin{document}

\maketitle

\renewcommand{\thefootnote}{\fnsymbol{footnote}}

\footnotetext[1]{CMAP Ecole Polytechnique, Route de Saclay, 91128 Palaiseau, 
France.}
\footnotetext[2]{DRI Institut Polytechnique des Sciences Avanc\'ees, 15-21 rue M. 
Grandcoing, 94200 Ivry-sur-Seine, France.}
\footnotetext[3]{Laboratoire J.A. Dieudonn\'e UMR 7351 CNRS-Universit\'e de 
Nice-Sophia Antipolis, Parc Valrose, 06108 Nice Cedex 02, France. 
\vspace{5pt} }

\footnotetext{
\textit{
M. Giacomini and O. Pantz are members of the DeFI team at INRIA 
Saclay \^Ile-de-France. \\
This work has been partially supported by the FMJH (Fondation 
Math\'ematique Jacques Hadamard) through a PGMO Gaspard Monge Project.}
\vspace{5pt} }

\footnotetext{\textit{e-mail:} \texttt{matteo.giacomini@polytechnique.edu; 
olivier.pantz@polytechnique.org; karim.trabelsi@ipsa.fr}}

\renewcommand{\thefootnote}{\arabic{footnote}}

\begin{abstract}
In this paper we introduce a novel certified shape optimization strategy - named 
Certified Descent Algorithm (CDA) - to account for the numerical error 
introduced by the Finite Element approximation of the shape gradient.
We present a goal-oriented procedure to derive a certified upper bound of the 
error in the shape gradient and we construct a fully-computable, constant-free 
\emph{a posteriori} error estimator inspired by the complementary energy 
principle.
The resulting CDA is able to identify a genuine descent direction at each 
iteration and features a reliable stopping criterion.
After validating the error estimator, some numerical simulations of the 
resulting certified shape optimization strategy are presented for the well-known 
inverse identification problem of Electrical Impedance Tomography.
\end{abstract}
%
%

\section{Introduction}
\label{ref:intro}

Shape optimization is a class of optimization problems in which the objective 
functional depends on the shape of the domain in which a Partial Differential 
Equation (PDE) is formulated and on the solution of the PDE itself.
Thus, we may view these problems as PDE-constrained optimization problems of a 
shape-dependent functional, the domain being the optimization variable and the 
PDE being the constraint. 
This class of problems has been tackled in the literature using both 
gradient-based and gradient-free methods and in this work we consider a strategy 
issue of the former group by computing the so-called shape gradient.

In most applications, the differential form of the objective functional with 
respect to the shape depends on the solution of a PDE which usually can only be 
solved approximately by means of a discretization strategy like the Finite 
Element Method. 
The approximation of the governing equation for the phenomenon under analysis 
introduces an uncertainty which may prevent the shape gradient from being 
strictly negative along the identified descent direction, that is, the 
approximated direction may not lead to any improvement of the objective 
functional we are trying to optimize. 
Moreover, due to the aforementioned approximation, stopping criteria based on 
the norm of the shape gradient may never be fulfilled if the \emph{a priori} 
given tolerance is too small with respect to the chosen discretization.
Within this framework, \emph{a posteriori} error estimators provide useful 
information to improve gradient-based algorithms for shape optimization.

Several works in the literature have highlighted the great potential of coupling 
\emph{a posteriori} error estimators to shape optimization algorithms.
In the pioneering work \cite{MeshRef-Banichuk}, the authors identify two 
different sources for the numerical error: on the one hand, the error arising 
from the approximation of the differential problem and on the other hand, the 
error due to the approximation of the geometry.
Starting from this observation, Banichuk \etal \ present a first attempt to use 
the information on the discretization of the differential problem provided by a 
recovery-based estimator and the error arising from the approximation of the 
geometry to develop an adaptive shape optimization strategy.
This work has been later extended by Morin \etal \ in \cite{COV:8787109}, where 
the adaptive discretization of the governing equations by means of the Adaptive 
Finite Element Method is linked to an adaptive strategy for the approximation of 
the geometry.
The authors derive estimators of the numerical error that are later used to 
drive an Adaptive Sequential Quadratic Programming algorithm to appropriately 
refine and coarsen the computational mesh.
Several other authors have used adaptive techniques for the approximation of 
PDE's in order to improve the accuracy of the solution and obtain better final 
configurations in optimal structural design problems. We refer to 
\cite{AFEM-Kikuchi, AFEM-Maute, Alauzet, Roche} for some examples. We remark 
that in all these works, \emph{a posteriori} estimators only provide qualitative 
information about the numerical error due to the discretization of the problems 
and are essentially used to drive mesh adaptation procedures.
To the best of our knowledge, no guaranteed fully-computable estimate has been 
investigated and the error in the shape gradient itself is not accounted for, 
thus preventing reliable stopping criteria to be derived.

In the present paper, the derivation of a fully-computable upper bound for the 
error in the shape gradient is tackled.
We neglect the contribution of the approximation of the geometry and we focus on 
the error arising from the discretization of the governing equation. 
The quantitative estimate of the error due to the numerical approximation of the 
shape gradient allows to identify a genuine descent direction and to introduce a 
reliable stopping criterion for the overall optimization strategy.
We propose a novel shape optimization strategy - named Certified Descent 
Algorithm (CDA) - that generates a sequence of minimizing shapes by certifying 
at each iteration the descent direction to be genuine and automatically stops 
when a reliable stopping criterion is fulfilled.

The rest of the paper is organized as follows. 
First, we introduce the general framework of shape optimization and shape 
identification problems and the so-called Boundary Variation Algorithm (Section 
\ref{ref:ShapePB}). 
In section \ref{ref:algorithm}, we account for the discretization error in the 
Quantity of Interest using the framework presented in \cite{Oden2001735}.
Then, we introduce the resulting Certified Descent Algorithm that couples the 
\emph{a posteriori} error estimator and the Boundary Variation Algorithm to 
derive a genuine descent direction for the shape optimization problem (Section 
\ref{ref:CDA}).
In section \ref{ref:InversePB}, we present the application of the Certified 
Descent Algorithm to the inverse problem of Electrical Impedance Tomography 
(EIT): after introducing the formulation of the identification problem as a 
shape optimization problem, we derive a fully-computable upper bound of the 
error in the shape gradient using the complementary energy principle.
Eventually, in section \ref{ref:numerics} we present some numerical tests of the 
application of CDA to the EIT problem and section \ref{ref:conclusion} 
summarizes our results.

\section{Shape optimization and shape identification problems}
\label{ref:ShapePB}

We consider an open domain $\Omega \subset \mathbb{R}^d$ ($d \geq 2$) with 
Lipschitz boundary $\partial\Omega$.
Let $V_\Omega$ be a separable Hilbert space depending on $\Omega$, we define 
$u_\Omega \in V_\Omega$ to be the solution of a state equation which is a 
linear PDE in the domain $\Omega$:
\begin{equation}
a_\Omega(u_\Omega,\delta u) = F_\Omega(\delta u) 
\qquad \forall \delta u \in V_\Omega
\label{eq:state}
\end{equation}
where $a_\Omega(\cdot,\cdot):V_\Omega \times V_\Omega \rightarrow \mathbb{R}$ 
is a continuous bilinear form satisfying the inf-sup condition
$$
\adjustlimits\inf_{w \in V_\Omega} \sup_{v \in V_\Omega} 
\frac{a_\Omega(v,w)}{\|v\| \|w\|} 
= \adjustlimits\inf_{v \in V_\Omega} \sup_{w \in V_\Omega} 
\frac{a_\Omega(v,w)}{\|v\| \|w\|} > 0 
$$
and $F_\Omega(\cdot)$ is a continuous linear form on $V_\Omega$, both of them 
depending on $\Omega$.
Under these assumptions, problem (\ref{eq:state}) has a unique solution 
$u_\Omega$.

We introduce a cost functional $J(\Omega)=j(\Omega,u_\Omega)$ which depends on 
the domain $\Omega$ itself and on the solution $u_\Omega$ of the state 
equation. We consider the following shape optimization problem
\begin{equation}
\min_{\Omega \in \mathcal{U}_{\text{ad}}} J(\Omega)
\label{eq:shapeOpt} 
\end{equation}
where $\mathcal{U}_{\text{ad}}$ is the set of admissible domains in 
$\mathbb{R}^d$.
Within this framework, problem (\ref{eq:shapeOpt}) may be viewed as a 
PDE-constrained optimization problem, in which we aim at minimizing the 
functional 
$j(\Omega,u)$ under the constraint $u=u_\Omega$, that is the minimizer $u$ is 
solution of the state equation (\ref{eq:state}). \\
In the following subsections, we recall the notion of shape gradient of 
$J(\Omega)$ in the direction $\theta$ and we apply the Steepest Descent Method 
to the shape optimization problem (\ref{eq:shapeOpt}).

\subsection{Differentiation with respect to the shape}

Let $X \subset W^{1,\infty}(\Omega;\mathbb{R}^d)$ be a Banach space and $\theta 
\in X$ be an admissible smooth deformation of $\Omega$. 
The cost functional $J(\Omega)$ is said to be $X$-differentiable at $\Omega \in 
\mathcal{U}_{\text{ad}}$ if there exists a continuous linear form $dJ(\Omega)$ 
on $X$ such that $\forall \theta \in X$
$$
J((\Id+\theta)\Omega) = J(\Omega) + \langle dJ(\Omega) , \theta \rangle + 
o(\theta).
$$
Several approaches are feasible to compute the shape gradient. Here we briefly 
recall the fast derivation method by C\'ea \cite{Cea1986} and the material 
derivative approach \cite{sokolowski1992introduction}.
Let us introduce the Lagrangian functional, defined for every admissible open 
set $\Omega$ and every $u, \ p \in V_\Omega$ by
\begin{equation}
\mathcal{L}(\Omega,u,p) = j(\Omega,u) + a_\Omega(u,p) - F_\Omega(p).
\label{eq:Lagrangian}
\end{equation}
Let $p_\Omega \in V_\Omega$ be the solution of the so-called adjoint problem, 
that is 
\begin{equation}
a_\Omega(\delta p,p_\Omega) + \left\langle \frac{\partial j}{\partial u}
(\Omega,u_\Omega), \delta p \right\rangle = 0 \qquad \forall \delta p \in 
V_\Omega .
\label{eq:adjoint}
\end{equation}
By applying the fast derivation method by C\'ea, we get the following expression 
of the shape gradient:
\begin{equation}
\langle dJ(\Omega),\theta \rangle = \left\langle 
\frac{\partial \mathcal{L}}{\partial\Omega}(\Omega,u_\Omega,p_\Omega), 
\theta \right\rangle.  
\label{eq:shapeDerCea} 
\end{equation}

An alternative procedure to compute the shape gradient relies on the definition 
of a diffeomorphism $\varphi: \mathbb{R}^d \rightarrow \mathbb{R}^d$ such that 
for every admissible set $\Omega$
$$
\Omega_{\varphi} \coloneqq \varphi(\Omega).
$$
Moreover, all functions $u,\ p \in V_\Omega$ defined on the reference domain 
$\Omega$ may be mapped to the deformed domain $\Omega_{\varphi}$ by
$$
u_\varphi \coloneqq u\circ\varphi^{-1} \qquad \text{and} \qquad 
p_\varphi\coloneqq p\circ\varphi^{-1}.
$$
We admit that $u \mapsto u_\varphi$ is a one-to-one map between $V_\Omega$ and 
$V_{\varphi(\Omega)}$.
The Lagrangian (\ref{eq:Lagrangian}) is said to admit a material derivative if 
there exists a linear form $\frac{\partial \mathcal{L}}{\partial \varphi}$ such 
that
$$
\mathcal{L}(\Omega_{\varphi},u_{\varphi},p_{\varphi}) = \mathcal{L}(\Omega,u,p) 
+ \left\langle \frac{\partial \mathcal{L}}{\partial \varphi}(\Omega,u,p) , 
\theta \right\rangle + o(\theta)
$$
where $\varphi=\Id+\theta$.
Provided that $u_\varphi$ is differentiable with respect to $\varphi$ at 
$\varphi=\Id$ in $V_{\varphi(\Omega)}$, from the fast derivation method of 
C\'ea we obtain the following expression for the shape gradient:
\begin{equation}
\langle dJ(\Omega),\theta \rangle = \left\langle 
\frac{\partial \mathcal{L}}{\partial \varphi}(\Omega,u_\Omega,p_\Omega) , 
\theta \right\rangle.
\label{eq:shapeDerDiff} 
\end{equation}

A variant of the latter method consists in computing the shape gradient via the 
Lagrangian formulation without explicitly constructing the material derivative 
of the state and adjoint solutions. We refer to \cite{LaurainSturm2015} for 
additional information about this approach.

\subsubsection*{Volumetric and surface expressions of the shape gradient}

The most common approach in the literature to compute the shape gradient is 
based on an Eulerian point of view and leads to a surface expression of the 
shape gradient. \\
The main advantage of this method relies on the fact that the boundary 
representation intuitively provides an explicit expression for the descent 
direction.
Let us assume that the shape gradient has the following form
$$
\langle dJ(\Omega),\theta \rangle = \int_{\partial\Omega}{h \theta \cdot n ds}
$$
then $\theta = - h n$ \ on \ $\partial\Omega$ is a descent direction.
Moreover, by Hadamard-Zol\'esio structure theorem it is well-known that the 
shape gradient is carried on the boundary of the shape 
\cite{Delfour:2001:SGA:501610} and using this approach the descent direction has 
to be defined only on $\partial\Omega$. 
Nevertheless, if the boundary datum of the state problem is not sufficiently 
smooth, the surface expression of the shape gradient may not exist or the 
descent direction $\theta$ may suffer from poor regularity. \\
Starting from the surface representation of the shape gradient, it is possible 
to derive a volumetric expression as well. 
Though the two expressions are equivalent in a continuous framework, they 
usually are not when considering their numerical counterparts, e.g. their Finite 
Element approximations: as a matter of fact, in \cite{HPS_bit} Hiptmair \etal \ 
prove that the volumetric formulation generally provides better accuracy when 
using the Finite Element Method.
Moreover, we may be able to derive the volumetric expression of the shape 
gradient even when its boundary representation fails to exist. \\
In this work, we consider the volumetric expression of the shape gradient in 
order to take advantage of the better accuracy it provides from a numerical 
point of view and to construct an estimator of the error in a Quantity of 
Interest using the procedure described by Oden and Prudhomme in 
\cite{Oden2001735}. \\
We remark that in order to derive a descent direction $\theta$ on $\Omega$ from 
the volumetric expression of the shape gradient, an additional variational 
problem has to be solved, as described in next subsection.

\subsection{The Boundary Variation Algorithm}
\label{ref:BVA}

From now on, we consider $X$ to be a Hilbert space. 
Starting from the formulation (\ref{eq:shapeDerDiff}), we seek a descent 
direction for the functional $J(\Omega)$. 
For this purpose, we solve an additional variational problem and we seek 
$\theta \in X$ such that
\begin{equation}
( \theta, \delta\theta )_X + \langle dJ(\Omega),\delta\theta \rangle = 0 \qquad 
\forall \delta\theta \in X.
\label{eq:variationalP}
\end{equation}

\begin{rmrk}
The choice of the scalar product in (\ref{eq:variationalP}) is a key point for 
the development of an efficient shape optimization method.
In \cite{tractionMethod}, the authors propose the so-called traction method to 
get rid of some irregularity issues in shape optimization problems. This 
approach is based on the regularization of the descent direction by means of a 
scalar product inspired by the linear elasticity equation.
In recent years, a comparison of the $L^2$, $H^1$ and $H^{-1}$ scalar products 
defined on a surface was presented \cite{Dogan20073898} 
but, as the authors state, the best choice is strongly dependent on the 
application of interest.
\end{rmrk}

In this section, we present the application of the Steepest Descent Method to a 
shape optimization problem.
After computing the solution of the state equation, we solve the adjoint problem 
to derive the expression of the shape gradient.
Then, a descent direction is identified through (\ref{eq:variationalP}) and is 
used to deform the domain. 
The resulting shape optimization strategy is known in the literature as Boundary 
Variation Algorithm \cite{smo-AP} and is sketched in script 
\ref{scpt:shape-opt-classic}.
\begin{lstlisting}[language=pseudo, escapeinside={/*@}{@*/}, 
caption={Continuous gradient method - The Boundary Variation Algorithm}, 
label=scpt:shape-opt-classic]
Given the domain $\Omega_0$, set $\ell=0$ and iterate:
1. Compute the solution $u_{\Omega_{\ell}}$ of the state equation;
2. Compute the solution $p_{\Omega_{\ell}}$ of the adjoint equation;
3. Compute a descent direction $\theta_{\ell} \in X$ solving $( \theta_{\ell}, 
\delta\theta )_X + \langle dJ(\Omega_{\ell}),\delta\theta \rangle = 0 \quad 
\forall \delta\theta \in X$;
4. Identify an admissible step $\mu_{\ell}$;
5. Update the domain $\Omega_{\ell+1} = (\Id + \mu_{\ell} 
\theta_{\ell})\Omega_{\ell}$;
6. While $| \langle dJ(\Omega_{\ell}),\theta_{\ell} \rangle | > 
\texttt{tol}$, $\ell=\ell+1$ and repeat.
\end{lstlisting}
We remark that this method relies on the computation at each iteration of a 
direction $\theta$ such that the shape gradient of the objective functional in 
this direction is strictly negative, that is we seek $\theta$ such that $\langle 
dJ(\Omega) , \theta \rangle < 0$.
In next subsection, we discuss the modifications that occur when moving from 
the continuous to the discretized formulation of the problems and consequently 
the conditions that the discretization of $\theta$ has to fulfill in order to be 
a genuine descent direction for the functional $J(\Omega)$.

\subsection{The discretized Boundary Variation Algorithm}

Let us denote by $u_\Omega^h$ and $p_\Omega^h$ the approximations of the 
state and adjoint equations arising from a Finite Element discretization. 
The discretized direction $\theta^h \in X$ is obtained by solving problem 
(\ref{eq:appTheta})
\begin{equation}
( \theta^h, \delta\theta )_X + \langle d_h J(\Omega),\delta\theta \rangle = 0 
\qquad \forall \delta\theta \in X
\label{eq:appTheta}
\end{equation}
where $\langle d_h J(\Omega),\delta\theta \rangle$ reads as follows:
\begin{equation}
\langle d_h J(\Omega),\delta\theta \rangle \coloneqq \left\langle \frac{\partial 
\mathcal{L}}{\partial \varphi}(\Omega,u_\Omega^h,p_\Omega^h), \delta\theta 
\right\rangle.
\label{eq:discreteShapeGrad}
\end{equation}
The discretized version of the Boundary Variation Algorithm is derived by 
substituting the continuous functions $u_\Omega$, $p_\Omega$ with their 
approximations $u_\Omega^h$, $p_\Omega^h$ and $\theta$ \ with \ $\theta^h$:

\begin{lstlisting}[language=pseudo, escapeinside={/*@}{@*/}, 
caption={Discretized gradient method - The discretized Boundary Variation 
Algorithm}, label=scpt:shape-opt-discrete]
Given the domain $\Omega_0$, set $\ell=0$ and iterate:
1. Compute the solution $u_{\Omega_{\ell}}^h$ of the state equation;
2. Compute the solution $p_{\Omega_{\ell}}^h$ of the adjoint equation;
3. Compute a descent direction $\theta_\ell^h \in X$ solving $( 
\theta_\ell^h, \delta\theta )_X + \langle d_h J(\Omega_{\ell}),\delta\theta 
\rangle = 0 \quad \forall \delta\theta \in X$;
4. Identify an admissible step $\mu_{\ell}$;
5. Update the domain $\Omega_{\ell+1} = (\Id + \mu_{\ell} 
\theta_\ell^h)\Omega_{\ell}$;
6. While $| \langle d_h J(\Omega_{\ell}),\theta_\ell^h \rangle | > 
\texttt{tol}$, $\ell=\ell+1$ and repeat.
\end{lstlisting}
We remark that due to the numerical error introduced by the Finite Element 
discretization, even though $\langle d_h J(\Omega),\theta^h \rangle < 0$, \ 
$\theta^h$ is not necessarily a genuine descent direction for the functional 
$J(\Omega)$. 
Moreover, it is important to notice that the stopping criterion (Algorithm 
\ref{scpt:shape-opt-discrete} - step 6) will usually not be fulfilled if the 
required tolerance $\texttt{tol}$ is too sharp with respect to the chosen 
discretization. \\
In order to bypass these issues, in the following sections we present a 
strategy to account for the error introduced by the approximation of the shape 
gradient. This results in a certification procedure that allows to verify 
whether a given direction is a genuine descent direction for the 
functional $J(\Omega)$ or not.

\subsection{Certification procedure for a genuine descent direction}
\label{ref:certification}

In this section we introduce the notion of certified descent direction, that is a 
direction which is verified to be a genuine descent direction for the functional 
$J(\Omega)$. 
As previously stated, a genuine descent direction $\theta$ is such that 
\begin{equation}
\langle d J(\Omega),\theta \rangle < 0 .
\label{eq:descentDirection}
\end{equation}
When moving from the continuous to the Finite Element framework, an approximation of the 
shape gradient is introduced (cf. definition (\ref{eq:discreteShapeGrad})) and consequently 
a numerical error appears. We define the error in the shape gradient $E^h$ as follows:
\begin{equation}
E^h \coloneqq \langle dJ(\Omega) - d_h J(\Omega),\theta^h \rangle.
\label{eq:errorH}
\end{equation}
By observing that $\langle d J(\Omega),\theta^h \rangle  = \langle d_h J(\Omega),\theta^h \rangle + E^h$, 
we get that a discretized direction $\theta^h$ is a descent direction for the objective 
functional $J(\Omega)$ if
\begin{equation}
\langle d_h J(\Omega),\theta^h \rangle + E^h < 0 .
\label{eq:condition-with-err}
\end{equation}
By construction, we seek $\theta^h$ such that $\langle d_h J(\Omega),\theta^h \rangle < 0$.
Nevertheless, this condition does not imply that $\theta^h$ is a genuine descent direction 
for $J(\Omega)$ since the quantity $E^h$ in (\ref{eq:condition-with-err}) may be either 
positive or negative.
In order to derive a relationship that stands independently from the sign of $E^h$ and 
since no \emph{a priori} information on the aforementioned sign is available, we modify 
(\ref{eq:condition-with-err}) by introducing the absolute value of $E^h$:
\begin{equation}
\langle d_h J(\Omega),\theta^h \rangle + E^h < \langle d_h J(\Omega),\theta^h \rangle + | E^h | < 0.
\label{eq:Grad+AbsErr}
\end{equation}
Let $\overline{E}$ be the upper bound of the quantity $| E^h |$.
From (\ref{eq:Grad+AbsErr}), we derive the following condition for the certification 
procedure: we seek $\theta^h$ such that
\begin{equation}
\langle d_h J(\Omega),\theta^h \rangle + \overline{E} < 0.
\label{eq:certified}
\end{equation}
It is straightforward to observe that if $\theta^h$ fulfills (\ref{eq:certified}), 
then it verifies (\ref{eq:descentDirection}) as well.
Thus, a direction $\theta^h$ fulfilling condition (\ref{eq:certified}) is said to be 
certified because it is guaranteed that the functional $J(\Omega)$ 
decreases along $\theta^h$.
We remark that the criterion (\ref{eq:certified}) ensures that $\theta^h$ is a 
genuine descent direction, whether it is the solution of equation 
(\ref{eq:appTheta}) or not. In particular, this also stands when the latter 
problem is only solved approximately. \\
In the following section, we present a strategy to construct a fully-computable 
guaranteed upper bound of the error in the approximation of the shape gradient 
in order to practically implement the certification procedure described above.

\section{Numerical error in the shape gradient}
\label{ref:algorithm}

In this section, we provide the detail of the technique used to derive an upper bound 
$\overline{E}$ of the error $| E^h |$ in the shape gradient.
The strategy to estimate the \emph{a posteriori} error 
in a Quantity of Interest (QoI) - namely, the shape gradient - is derived from the 
work \cite{Oden2001735} by Oden and Prudhomme. 
Basic idea relies on the definition of an adjoint problem whose right-hand side 
is the quantity whose error estimate is sought.

\subsection{Bound for the approximation error of a linear functional}

Here we briefly recall the aforementioned strategy by Oden and Prudhomme for the derivation of an 
\emph{a posteriori} error estimate for a bounded linear functional $Q: V_\Omega \rightarrow \mathbb{R}$, 
also known as Quantity of Interest.
Let $e_u \coloneqq u_\Omega-u_\Omega^h$ be the error between the function $u_\Omega$ solution of the 
state problem (\ref{eq:state}) and its Finite Element counterpart $u_\Omega^h$. 
We are interested in evaluating the target $Q(u_\Omega)$ and the accuracy of its approximation $Q(u_\Omega^h)$ 
is expressed via the following quantity:
$$
E_Q \coloneqq Q(u_\Omega) - Q(u_\Omega^h) = Q(u_\Omega-u_\Omega^h) = Q(e_u)
$$
where the first equality follows from the linearity of $Q$. 
In order to compute the error $Q(e_u)$, we introduce the residue associated with the approximation of the 
state problem (\ref{eq:state}):
\begin{equation}
\mathcal{R}_\Omega^u(\delta u) \coloneqq F_\Omega(\delta u) - a_\Omega(u_\Omega^h,\delta u).
\label{eq:residueGoalOriented}
\end{equation}
Moreover, we recall that the error $e_u$ is solution of the so-called residual equation that reads
\begin{equation}
a_\Omega(e_u,\delta u) = \mathcal{R}_\Omega^u(\delta u)
\qquad \forall \delta u \in V_\Omega .
\label{eq:residalEqGoalOriented}
\end{equation}
We highlight that (\ref{eq:residueGoalOriented}) contains all the information related to the numerical error $e_u$, thus 
the evaluation of $Q(e_u)$ reduces to the derivation of a relationship between $\mathcal{R}_\Omega^u$ 
and $Q(e_u)$ itself. Hence we seek a so-called influence function that carries the information on the 
effect of the residue - i.e. the effect of the error $e_u$ - on the functional $Q$. 
Let us assume that there 
exists an influence function $r_\Omega$ such that $Q(e_u) = \mathcal{R}_\Omega^u(r_\Omega)$.
We may identify $r_\Omega$ with its Riesz element and owing to (\ref{eq:residalEqGoalOriented}), we get 
$$
Q(e_u) = \mathcal{R}_\Omega^u(r_\Omega) = a_\Omega(e_u,r_\Omega) .
$$
By combining the above information we are able to derive the following Boundary Variation Problem - 
known as adjoint problem - to construct the influence function. 
In particular, we seek $r_\Omega \in V_\Omega$ such that
\begin{equation}
a_\Omega(\delta r,r_\Omega) = Q(\delta r)
\qquad \forall \delta r \in V_\Omega .
\label{eq:adjointGoalOriented}
\end{equation}
Existence and uniqueness of the solution $r_\Omega$ of the adjoint problem follow from the Lax-Milgram 
theorem. From a practical point of view, the solution of the aforementioned adjoint problem is 
only approximated - usually using the same method as for the state problem - and an additional error 
$e_r \coloneqq r_\Omega-r_\Omega^h$ is introduced. 
Hence, owing to the Galerkin orthogonality the evaluation of $Q(e_u)$ reads as
\begin{equation}
Q(e_u) = a_\Omega(e_u,r_\Omega) = a_\Omega(e_u,r_\Omega-r_\Omega^h) =  a_\Omega(e_u,e_r) .
\label{eq:ErrorTargetGoalOriented}
\end{equation}

\subsection{Variational formulation of the error in the shape gradient}
\label{ref:errorQoI}

Several works in the literature \cite{2006-nonlin, nonlin-goal} have dealt with an 
extension of the aforementioned framework to compute a bound of the approximation error 
in a target functional to the case of non-linear Quantities of Interest. 
To account for the error (\ref{eq:errorH}), we follow the approach proposed by these 
authors by performing a linearization of the functional whose error estimate is 
sought. 
Thus we rewrite the numerical error (\ref{eq:errorH}) in the shape gradient 
by introducing the linearized error $\widetilde{E}^h$:
\begin{equation}
\begin{aligned}
E^h = & \left\langle \frac{\partial \mathcal{L}}{\partial \varphi}
(\Omega,u_\Omega,p_\Omega) - \frac{\partial\mathcal{L}}{\partial\varphi}
(\Omega,u_\Omega^h,p_\Omega^h), \theta^h \right\rangle \\
 & \simeq \frac{\partial^2 \mathcal{L}}{\partial \varphi \partial u}
 (\Omega,u_\Omega^h,p_\Omega^h)[\theta^h,u_\Omega-u_\Omega^h] 
+ \frac{\partial^2 \mathcal{L}}{\partial \varphi \partial p}
(\Omega,u_\Omega^h,p_\Omega^h)[\theta^h,p_\Omega-p_\Omega^h] \eqqcolon 
\widetilde{E}^h.
\end{aligned}
\label{eq:errDJ}
\end{equation}

In order to compute an upper bound of the error (\ref{eq:errDJ}), 
we introduce two adjoint problems, each of which is associated 
with one term on the right-hand side of (\ref{eq:errDJ}).
Thus, we seek $r_\Omega,s_\Omega \in V_\Omega$ such that
\begin{equation}
\begin{aligned}  
& a_\Omega(\delta r,r_\Omega) = \frac{\partial^2 \mathcal{L}}
{\partial\varphi \partial u}(\Omega,u_\Omega^h,p_\Omega^h)[\theta^h,\delta 
r] \qquad \forall \delta r \in V_\Omega \\
& a_\Omega(\delta s,s_\Omega) = \frac{\partial^2 \mathcal{L}}
{\partial\varphi \partial p}(\Omega,u_\Omega^h,p_\Omega^h)[\theta^h,\delta 
s] \qquad \forall \delta s \in V_\Omega 
\end{aligned}
\label{eq:adjointDJ}
\end{equation}
We remark that in order for the aforementioned adjoint problems to be 
well-posed, their right-hand sides have to be linear and continuous forms on 
$V_\Omega$ and this motivates the linearization introduced in (\ref{eq:errDJ}).

Let us denote by $r_\Omega^h$, $s_\Omega^h$ the approximations of the solutions 
$r_\Omega$, $s_\Omega$ of equations (\ref{eq:adjointDJ}) arising from a Finite 
Element discretization.
By plugging (\ref{eq:adjointDJ}) into (\ref{eq:errDJ}), we may derive the 
following upper bound $\overline{E}$ for the numerical error in the shape 
gradient:
\begin{equation}
\begin{aligned}  
|E^h| \simeq |\widetilde{E}^h| & \leq 
|a_\Omega(u_\Omega-u_\Omega^h,r_\Omega-r_\Omega^h) | + 
|a_\Omega(p_\Omega-p_\Omega^h,s_\Omega-s_\Omega^h)| \\
& \leq \enorm{u_\Omega-u_\Omega^h}\enorm{r_\Omega-r_\Omega^h} + 
\enorm{p_\Omega-p_\Omega^h}\enorm{s_\Omega-s_\Omega^h}  \eqqcolon 
\overline{E}
\end{aligned}
\label{eq:upperBound}
\end{equation}
where $\enorm{\cdot}$ is the energy-norm induced by the bilinear form 
$a_\Omega(\cdot,\cdot)$.
The first inequality follows from triangle inequality and Galerkin 
orthogonality whereas the upper bound $\overline{E}$ is derived exploiting 
Cauchy-Schwarz inequality.

\begin{rmrk}
In (\ref{eq:upperBound}) we derived an upper bound for the numerical error in 
the linearized shape gradient and not in the shape gradient itself.
For the rest of this paper, we will follow the framework in \cite{2006-nonlin} 
by assuming the linearization error to be negligible and $\overline{E}$ to be 
an upper bound of the numerical error $E^h$ itself and not of its linearized 
version $\widetilde{E}^h$.
In section \ref{ref:validation}, a validation of the error estimator is 
presented for the case of Electrical Impedance Tomography: we will verify that 
the linearization error is indeed negligible with respect to the error due to 
the Finite Element discretization and thus the previous assumption stands.
\end{rmrk}

In order to fully compute the error estimator (\ref{eq:upperBound}), we 
have to estimate the energy-norm of the error for:
\begin{itemize}
 \item the state equation;
 \item the adjoint equation used to compute the shape gradient;
 \item the two adjoint equations associated with the Quantity of Interest.
\end{itemize}
Several strategies are possible to tackle these issues. In this paper, we 
propose a method inspired by the complementary energy principle which allows to 
derive fully-computable, constant-free estimators by solving an 
additional variational problem for each term under analysis in order to retrieve 
a flux estimate. 
These estimates are problem-dependent and will be detailed in section 
\ref{ref:InversePB} when presenting the case of Electrical Impedance Tomography.

\section{The Certified Descent Algorithm}
\label{ref:CDA}

We are now ready to introduce the novel Certified Descent Algorithm, arising 
from the coupling of the Boundary Variation Algorithm for shape optimization 
(Section \ref{ref:BVA}) and the goal-oriented estimator for the error in the 
shape gradient (Section \ref{ref:errorQoI}). 
In script \ref{scpt:shape-opt-adaptive}, we present a variant of algorithm 
\ref{scpt:shape-opt-classic} that takes advantage of the previously introduced 
\emph{a posteriori} estimator for the error in the shape gradient in order to 
bypass the issues due to the discretization of the problem.

First, the procedure constructs an initial computational domain. At each 
iteration, the algorithm solves the state and adjoint problems [steps 1 and 2] 
and computes a descent direction $\theta^h$ solving equation (\ref{eq:appTheta}) 
[step 3]. 
Then, the adjoint problems (\ref{eq:adjointDJ}) are solved and an upper bound of 
the numerical error in the shape gradient along the direction $\theta^h$ is 
computed [step 4]. If condition (\ref{eq:certified}) is not fulfilled, the mesh 
is adapted in order to improve the error estimate.
This procedure is iterated  until the direction $\theta^h$ is a certified 
descent direction for $J(\Omega)$ [step 5]. 
Once a certified descent direction has been identified, we compute a step 
[step 6] such that the following Armijo condition is fulfilled: let us consider 
the iteration $\ell$, given $0 < \alpha < 1$  we use a backtracking strategy to 
identify the step $\mu_{\ell} \in \mathbb{R}^+$ such that
$$
J \left( \left( \Id + \mu_{\ell} \theta_\ell^h \right)\Omega_{\ell} \right) 
\leq J \left(\Omega_{\ell} \right) + \alpha \mu_{\ell} \langle d_h J 
\left(\Omega_{\ell} \right),\theta_\ell^h \rangle.
$$
An alternative bisection-based line search technique has been proposed by Morin 
\etal \ in \cite{COV:8787109}. \\
Then the shape of the domain is updated according to the computed perturbation 
of the identity $\Id + \mu_{\ell} \theta_\ell^h$ [step 7]. 
Eventually, a novel stopping criterion is proposed [step 8] in order to use the 
information embedded in the error bound $\overline{E}$ to derive a reliable 
condition to end the evolution of the algorithm.

\begin{lstlisting}[language=pseudo, escapeinside={/*@}{@*/}, 
caption={Discretized gradient method - The Certified Descent Algorithm}, 
label=scpt:shape-opt-adaptive]
Given the domain $\Omega_0$, set $\ell=0$ and iterate:
1. Compute the solution $u_{\Omega_{\ell}}^h$ of the state equation;
2. Compute the solution $p_{\Omega_{\ell}}^h$ of the adjoint equation;
3. Compute a descent direction $\theta_\ell^h \in X$ solving $( 
\theta_\ell^h, \delta\theta )_X + \langle d_h J(\Omega_{\ell}),\delta\theta 
\rangle = 0 \quad \forall \delta\theta \in X$;
4. Compute an upper bound $\overline{E}$ of the numerical error $|E^h|$:
   (a) Compute the solutions $r_{\Omega_{\ell}}^h$ and $ s_{ \Omega_{ \ell 
}}^h$ to estimate the error in the QoI;
   (b) Compute $\overline{E} = \enorm{u_{\Omega_{\ell}}-u_{\Omega_{\ell}}^h}
   \enorm{r_{\Omega_{\ell}}-r_{\Omega_{\ell}}^h} 
+ \enorm{p_{\Omega_{\ell}}-p_{\Omega_{\ell}}^h}\enorm{s_{\Omega_{\ell}}-s_{\Omega_{\ell}}^h}$;
5. If $\langle d_hJ(\Omega_{\ell}),\theta_\ell^h \rangle + \overline{E} \geq  
0$, refine the mesh and go to 1;
6. Identify an admissible step $\mu_{\ell}$;
7. Update the shape $\Omega_{\ell+1} = (\Id + \mu_{\ell} \theta_\ell^h) 
\Omega_{\ell}$;
8. While $| \langle d_hJ(\Omega_{\ell}),\theta_\ell^h \rangle | + \overline{E} 
>  \texttt{tol}$, $\ell=\ell+1$ and repeat.
\end{lstlisting}
The advantage of the Certified Descent Algorithm is twofold.
On the one hand, the computation of the upper bound of the numerical error 
in the shape gradient provides useful information to identify a certified 
descent direction at each iteration of the optimization algorithm and to 
construct a certified shape optimization strategy. 
On the other hand, the fully-computable and constant-free error estimator provides  
quantitative information to derive a reliable stopping criterion 
for the overall optimization procedure.

The novelty of algorithm \ref{scpt:shape-opt-adaptive} is the certification procedure that plays 
a crucial role in steps 4, 5 and 8. 
A key aspect of the procedure is the mesh adaptation routine that has to be run if condition 
(\ref{eq:certified}) is not fulfilled by the current configuration.
Owing to the subsequent refinements of the mesh (cf. Algorithm \ref{scpt:shape-opt-adaptive} - step 5), 
the approximated solution 
$u_\Omega^h$ tends to the exact solution $u_\Omega$ and analogously do the solutions $p_\Omega^h$, 
$r_\Omega^h$ and $s_\Omega^h$ of the adjoint problems. 
Thus the term $\overline{E}$ tends to zero as 
the mesh size tends to zero, assuring that (\ref{eq:certified}) is eventually fulfilled. 
From a practical point of view, in order to guarantee that condition (\ref{eq:certified}) is 
fulfilled in a reasonable number of iterations of the adaptation routine, we construct a 
refinement strategy to explicitly reduce the error $\overline{E}$ at each iteration.
In particular, we perform goal-oriented mesh adaptation as suggested in \cite{Oden2001735}: 
at each iteration, we construct the upper bound $\overline{E}$ and an indicator 
based on the estimator of the error in the shape gradient.
This approach exploits the previously constructed estimator 
to localize the areas of the domain that are mainly responsible for the error in the Quantity of 
Interest and performs a targeted refinement in order to concurrently reduce the error in the shape 
gradient and limit the number of newly inserted Degrees of Freedom.
The efficiency of this strategy has been extensively studied in the literature 
\cite{MR1957989, MR2101976} and the results in section \ref{ref:numerics} confirm 
the ability of this method to reduce the targeted error.

\section{An inverse identification problem: the case of Electrical Impedance 
Tomography}
\label{ref:InversePB}

We present the application of the Certified Descent Algorithm 
to the problem of Electrical Impedance Tomography.
The choice of EIT as test case has to be interpreted as a proof 
of concept to preliminarily assess the validity of the discussed 
method on a non-trivial scalar problem before studying the vectorial case. 

Let us consider an open domain $\mathcal{D} \subset \mathbb{R}^2$. We suppose 
that there exists an open subdomain $\Omega \subset\subset \mathcal{D}$ such 
that some given physical properties of the problem under analysis are 
discontinuous along the interface $\partial\Omega$ between the inclusion 
$\Omega$ and the complementary set $\mathcal{D} \setminus \Omega$. 
The location and the shape of the inclusion are to be determined, thus $\Omega$ 
acts as unknown parameter in the state equations and in the inversion procedure. 
Our aim is to identify the inclusion $\Omega$ by performing non-invasive 
measurements on the boundary $\partial\mathcal{D}$ of the domain $\mathcal{D}$. 
This problem is well-known in the literature and is often referred to as 
Calder\'on's problem. Several review papers on Electrical Impedance Tomography 
have been published in the literature over the years. We refer to 
\cite{Calderon1980, Cheney99electricalimpedance, 0266-5611-18-6-201} for more 
details on the physical problem, its mathematical formulation and its numerical 
approximation. \\
Let $\chi_\Omega$ be the characteristic function of the open set $\Omega$, we 
define the conductivity $k_\Omega$ as a piecewise constant function such that 
$k_\Omega \coloneqq k_I \chi_\Omega + k_E (1-\chi_\Omega), \ k_I, \ k_E 
>0$.
We introduce two Boundary Value Problems on the domain $\mathcal{D}$, 
respectively with Neumann and Dirichlet boundary conditions on 
$\partial\mathcal{D}$:
\begin{align}
\left\{
\begin{aligned}
& -k_\Omega \Delta u_{\Omega,N} + u_{\Omega,N} = 0 \quad & \text{in} \ 
\mathcal{D} \setminus \partial\Omega \\
& \llbracket u_{\Omega,N} \rrbracket = 0 \quad & \text{on} \ \partial\Omega \\
& \llbracket k_\Omega \nabla u_{\Omega,N} \cdot n \rrbracket = 0 \quad & 
\text{on} \ \partial\Omega \\
& k_E \nabla u_{\Omega,N} \cdot n = g \quad & \text{on} \ \partial\mathcal{D}
\end{aligned}
\right.
\label{eq:Neumann}
\end{align}
\begin{align}
\left\{
\begin{aligned}
& -k_\Omega \Delta u_{\Omega,D} + u_{\Omega,D} = 0 \quad & \text{in} \ 
\mathcal{D} \setminus \partial\Omega \\
& \llbracket u_{\Omega,D} \rrbracket = 0 \quad & \text{on} \ \partial\Omega \\
& \llbracket k_\Omega \nabla u_{\Omega,D} \cdot n \rrbracket = 0 \quad & 
\text{on} \ \partial\Omega \\
& u_{\Omega,D} = U_D \quad & \text{on} \ \partial\mathcal{D}
\end{aligned}
\right.
\label{eq:Dirichlet}
\end{align}
where the boundary data $g \in L^2(\partial\mathcal{D})$ and $U_D \in 
H^{\frac{1}{2}}(\partial\mathcal{D})$ arise from the performed physical 
measurements. 
As previously stated, we are interested in identifying the shape and the 
location of the inclusion, fitting given boundary measurements $g$ and $U_D$ of 
the flux and the potential.

\subsection{State problems}
\label{ref:state}

In order to approximate problems (\ref{eq:Neumann}) and (\ref{eq:Dirichlet}) by 
means of the Finite Element Method, first we introduce their variational 
formulations. \\
Let $a_\Omega(\cdot,\cdot)$ be the bilinear form associated with both the 
problems and $F_{\Omega,i}(\cdot), \ i=N,D$ the linear forms respectively for 
the Neumann and the Dirichlet problem:
\begin{gather}
a_\Omega(u,\delta u) = \int_{\mathcal{D}}{\Big(k_\Omega \nabla u \cdot 
\nabla \delta u + u \delta u \Big) dx} 
\label{eq:a} \\
F_{\Omega,N}(\delta u) = \int_{\partial\mathcal{D}}{g \delta u \ ds} \qquad 
\text{and} \qquad F_{\Omega,D}(\delta u) = 0.
\label{eq:F}
\end{gather}
We consider $u_{\Omega,N},\ u_{\Omega,D} \in H^1(\mathcal{D})$ such that 
$u_{\Omega,D}=U_D \ \text{on} \ \partial \mathcal{D}$, solutions of the 
following Neumann and Dirichlet variational problems $\forall \delta u_N \in 
H^1(\mathcal{D})$ and $\forall \delta u_D \in H^1_0(\mathcal{D})$:
\begin{equation}
a_\Omega(u_{\Omega,i},\delta u_i) = F_{\Omega,i}(\delta u_i) 
\quad , \quad i=N,D.
\label{eq:stateEIT}
\end{equation}
For the Dirichlet problem, the non-homogeneous boundary datum is taken care of 
by means of a classical substitution technique.
The corresponding discretized formulations of (\ref{eq:stateEIT}) may be derived 
by replacing the analytical solutions $u_{\Omega,N}$ and $u_{\Omega,D}$ with 
their approximations $u_{\Omega,N}^h$ and $u_{\Omega,D}^h$ which belong to the 
space of Lagrangian Finite Element functions.
In a similar fashion, $\theta^h$ is the solution of equation (\ref{eq:appTheta}) 
computed using a Lagrangian Finite Element space.
The degree chosen for the Finite Element basis functions will be discussed in 
section \ref{ref:numerics}.

\subsection{A shape optimization approach}
\label{ref:shapeOpt}

Let us consider the Kohn-Vogelius functional first introduced in 
\cite{Wexler:85} and later investigated by Kohn and Vogelius in 
\cite{CPA:CPA3160400605}:
\begin{equation}
J(\Omega) = \frac{1}{2}\int_{\mathcal{D}}{\Big( k_\Omega \left| 
\nabla(u_{\Omega,N} - u_{\Omega,D})\right|^2 + |u_{\Omega,N} - u_{\Omega,D}|^2  
\Big) dx}.
\label{eq:kohn-vogelius}
\end{equation}
In (\ref{eq:kohn-vogelius}), $u_{\Omega,N}$ and $u_{\Omega,D}$ 
respectively stand  for the solutions of the state problems (\ref{eq:Neumann}) 
and (\ref{eq:Dirichlet}). 
Owing to (\ref{eq:a}), we may rewrite the objective functional 
(\ref{eq:kohn-vogelius}) as:
$$
J(\Omega) = \frac{1}{2} a_\Omega\Big( u_{\Omega,N} - u_{\Omega,D} , 
u_{\Omega,N} - u_{\Omega,D} \Big).
$$
The inverse identification problem of Electrical Impedance Tomography may be 
written as the PDE-constrained optimization problem (\ref{eq:shapeOpt}) in 
which we seek the open subset $\Omega$ that minimizes (\ref{eq:kohn-vogelius}). 
\\
In order to solve this problem, we consider the Certified Descent 
Algorithm using the shape gradient of $J(\Omega)$.
First, we need to determine the adjoint solutions $p_{\Omega,N}$ and 
$p_{\Omega,D}$ associated with the states $u_{\Omega,N}$ and $u_{\Omega,D}$: 
the Kohn-Vogelius problem is self-adjoint and we get that 
$p_{\Omega,N}=u_{\Omega,N} - u_{\Omega,D}$ and $p_{\Omega,D}=0$.
Let $\theta \in W^{1,\infty}(\mathcal{D};\mathbb{R}^2)$ be an admissible 
deformation of the domain such that $\theta=0 \ \text{on} \ 
\partial\mathcal{D}$. 
As previously mentioned, the most common approach in the literature to compute 
the shape gradient leads to the surface expression
\begin{equation}
\langle dJ(\Omega),\theta \rangle = \frac{1}{2}\int_{\partial\Omega}{ \Big(  
\llbracket k_\Omega \rrbracket \Big( \Big| \frac{\partial 
u_{\Omega,N}}{\partial \tau} \Big|^2 - \Big| \frac{\partial 
u_{\Omega,D}}{\partial \tau} \Big|^2 \Big) 
- \llbracket k_\Omega^{-1} \rrbracket \Big( \Big| k_\Omega \frac{\partial  
u_{\Omega,N}}{\partial n} \Big|^2 - \Big| k_\Omega \frac{\partial 
u_{\Omega,D}}{\partial n} \Big|^2 \Big) \Big)(\theta \cdot n) ds}
\label{eq:surfaceKV}
\end{equation}
where $n$ is the outward normal to $\partial\Omega$, $\tau$ is the tangential 
direction to $\partial\Omega$ and $\llbracket k_\Omega \rrbracket = k_E -k_I$ 
and $\llbracket k_\Omega^{-1} \rrbracket = k_E^{-1} -k_I^{-1}$ are the jumps 
across $\partial\Omega$. 
Let us now introduce the following operator:
\begin{equation}
\langle G(\Omega,u) , \theta \rangle = \frac{1}{2}\int_{\mathcal{D}}{  
\Big(k_\Omega M(\theta) \nabla u \cdot \nabla u - \nabla \cdot \theta \ u^2 
\Big) dx} 
\label{eq:operatorG}
\end{equation}
where $M(\theta) = \nabla \theta + \nabla \theta^T - (\nabla \cdot \theta) \Id$. 
From (\ref{eq:shapeDerDiff}), we get the volumetric expression of the shape  
gradient of (\ref{eq:kohn-vogelius}):
\begin{equation}
\langle dJ(\Omega),\theta \rangle = \langle G(\Omega,u_{\Omega,N}) -  
G(\Omega,u_{\Omega,D}) , \theta \rangle.
\label{eq:volumetricKV}
\end{equation}
We refer to \cite{PANTZ-sauts} for more details on the differentiation of the 
Kohn-Vogelius functional and its application to the identification of 
discontinuities of the conductivity parameter.

\subsection{\emph{A posteriori} error estimate for the shape gradient}

Since the Kohn-Vogelius problem is self-adjoint, (\ref{eq:upperBound}) reduces 
to
\begin{equation}
\overline{E} = \enorm{u_{\Omega,N}- u_{\Omega,N}^h}\enorm{r_{\Omega,N}-r_{\Omega,N}^h} + 
\enorm{u_{\Omega,D}-u_{\Omega,D}^h}\enorm{r_{\Omega,D}-r_{\Omega,D}^h}
\label{eq:UBeit} 
\end{equation}
where $r_{\Omega,N}$ and $r_{\Omega,D}$ are the solutions of the adjoint 
problems introduced to evaluate the contributions of the Neumann and Dirichlet 
state problems to the error in the Quantity of Interest. Thus, we seek 
$r_{\Omega,N} \in H^1(\mathcal{D})$ and $r_{\Omega,D} \in H^1_0(\mathcal{D})$ 
such that respectively $\forall \delta r_N \in H^1(\mathcal{D})$ and $\forall 
\delta r_D \in H^1_0(\mathcal{D})$
\begin{equation}
a_\Omega(\delta r_i,r_{\Omega,i}) = H_{\Omega,i}(\delta r_i) 
\quad , \quad i=N,D
\label{eq:adjointEIT}
\end{equation}
where the linear forms on the right-hand sides of equations 
(\ref{eq:adjointEIT}) read as
\begin{equation}
H_{\Omega,i}(\delta r) = \frac{\partial G}{\partial u}(\Omega,u_{\Omega,i}^h)[ 
\theta^h,\delta r] \quad , \quad i=N,D.
\label{eq:Fadj}
\end{equation}

In order to obtain a computable upper bound for the error in the shape gradient, 
we seek an estimate of the energy-norm of the error for the state and adjoint 
solutions in (\ref{eq:UBeit}). \\
Let $e_{\Omega,i} = u_{\Omega,i} - u_{\Omega,i}^h$ and $\epsilon_{\Omega,i} = 
r_{\Omega,i} - r_{\Omega,i}^h$ for $i=N,D$. 
In the following subsections, we derive the estimates of the energy-norm of the 
$e_{\Omega,i}$'s and the $\epsilon_{\Omega,i}$'s using a strategy inspired by 
the so-called complementary energy principle 
\cite{Repin:2000:PEE:343641.343640}. 
In practice, we introduce a dual flux variable for every problem and each bound 
is computed by solving an additional adjoint problem thus leading to a better 
approximation of the numerical fluxes $\nabla e_{\Omega,i}$'s and $\nabla 
\epsilon_{\Omega,i}$'s.
For additional information on this approach, we refer to 
\cite{Vejchodsky20137194}.

\subsubsection{Error estimates based on the complementary energy principle: the 
case of the state equations}
\label{ref:stateEstimates}

For the sake of readability, let us rename $H^1(\mathcal{D})$ as $V_N$ and  
$H^1_0(\mathcal{D})$ as $V_D$.
We recall the previously mentioned residual equations such that $\forall \delta u_N \in V_N$ and  
$\forall \delta u_D \in V_D$
\begin{equation}
a_\Omega(e_{\Omega,i},\delta u_i) = F_{\Omega,i}(\delta u_i) - a_\Omega 
(u_{\Omega,i}^h,\delta u_i) \quad , \quad i=N,D.
\label{eq:residualEq} 
\end{equation}
We recall that solving equation (\ref{eq:residualEq}) is equivalent to the  
following minimization problem, that is we seek $w \in V_i$ such that
\begin{equation}
E_{\Omega,i}(e_{\Omega,i}) = \min_{w \in V_i} E_{\Omega,i}(w) 
\quad , \quad i=N,D
\label{eq:minW}
\end{equation}
where the global energy functional associated with the Neumann and Dirichlet 
problems reads
\begin{equation}
E_{\Omega,i}(w) = \frac{1}{2} \int_{\mathcal{D}}{\Big(k_\Omega | \nabla w |^2 
+ |w|^2 \Big) dx} + \int_{\mathcal{D}}{\Big(k_\Omega \nabla u_{\Omega,i}^h 
\cdot \nabla w + u_{\Omega,i}^h w \Big) dx} - F_{\Omega,i}(w).
\label{eq:EnergyMech}
\end{equation}
By introducing an additional variable $z = \nabla w$ and a dual variable 
$\sigma_{\Omega,i} \in L^2(\mathcal{D};\mathbb{R}^d)$, we may construct the 
Lagrangian functional $L_{\Omega,i}:V_i \times L^2(\mathcal{D}) \times 
L^2(\mathcal{D};\mathbb{R}^d) \rightarrow \mathbb{R}$ which has the following 
form
\begin{equation}
\begin{aligned}
L_{\Omega,i}(w,z,\sigma_{\Omega,i}) = & \frac{1}{2} \int_{\mathcal{D}}{\Big( 
k_\Omega |z|^2 + |w|^2 \Big) dx} + \int_{\mathcal{D}}{\Big(k_\Omega \nabla 
u_{\Omega,i}^h \cdot z + u_{\Omega,i}^h w \Big) dx} - F_{\Omega,i}(w) \\ 
& + \int_{\mathcal{D}}{\sigma_{\Omega,i} \cdot (\nabla w - z) dx}.
\end{aligned}
\label{eq:lagrangianEIT}
\end{equation}
Thus the minimization problem (\ref{eq:minW}) may be rewritten as a min-max  
problem and owing to the Lagrange duality, we get
\begin{equation}
\begin{aligned}
 \min_{w \in V_i} E_{\Omega,i}(w) & = \min_{\mytop{w \in V_i}{z = \nabla w}}  
\max_{\sigma_{\Omega,i} \in L^2(\mathcal{D};\mathbb{R}^d)} 
L_{\Omega,i}(w,z,\sigma_{\Omega,i}) \\
 & = \max_{\sigma_{\Omega,i} \in L^2(\mathcal{D};\mathbb{R}^d)} \min_{\mytop{w 
\in V_i}{z =\nabla w}} L_{\Omega,i}(w,z,\sigma_{\Omega,i}).
\end{aligned}
\label{eq:minMax}
\end{equation}
We consider the space $\Hdiv=\{ \sigma \in L^2(\mathcal{D}; 
\mathbb{R}^d) \ : \ \nabla \cdot \sigma \in L^2(\mathcal{D})\}$. Let 
$\sigma_{\Omega,i} \in \Hdiv$ for $i=N,D$, from the system of first-order 
optimality conditions for $L_{\Omega,i} (w,z,\sigma_{\Omega,i})$ we derive the 
following relationships among the variables:
\begin{align}
\left\{
\begin{aligned}
& z = k_\Omega^{-1}\sigma_{\Omega,i} - \nabla u_{\Omega,i}^h \quad & \text{in} 
\ \mathcal{D} \\
& w = \nabla \cdot \sigma_{\Omega,i} - u_{\Omega,i}^h \quad & \text{in} \ 
\mathcal{D} \\
& \sigma_{\Omega,N} \cdot n = g  \quad & \text{on} \ \partial\mathcal{D} 
\end{aligned}
\right.
\label{eq:Optimality}
\end{align}
Hence, by plugging (\ref{eq:Optimality}) into (\ref{eq:minMax}), we get the 
following maximization problems for $i=N,D$
\begin{equation}
E_{\Omega,i}(e_{\Omega,i}) = \max_{\mytop{\sigma_{\Omega,i} \in 
\Hdiv}{(\sigma_{\Omega,N} \cdot n = g )}} -\frac{1}{2} 
\int_{\mathcal{D}}{\Big(k_\Omega^{-1} | \sigma_{\Omega,i} - k_\Omega \nabla 
u_{\Omega,i}^h|^2 + |\nabla \cdot \sigma_{\Omega,i} - u_{\Omega,i}^h|^2 \Big) 
dx} 
\label{eq:maxS}
\end{equation}
where the objective functional is known as the dual complementary energy 
associated with the problems.

In order to compute the dual flux variables, we derive the first-order 
optimality conditions for the dual complementary energy functional in 
(\ref{eq:maxS}).
Thus, we seek $\sigma_{\Omega,N},\ \sigma_{\Omega,D} \in 
\Hdiv$ such that $\sigma_{\Omega,N} \cdot n = g \ \text{on} \ \partial 
\mathcal{D}$ which satisfy $\forall \delta\sigma_N,\ \delta\sigma_D \in \Hdiv$ 
such that $\delta\sigma_N \cdot n=0  \ \text{on} \ \partial \mathcal{D}$
\begin{equation}
\int_{\mathcal{D}}{\Big(k_\Omega^{-1} \sigma_{\Omega,i}\cdot\delta\sigma_i + 
(\nabla \cdot \sigma_{\Omega,i} )( \nabla \cdot \delta\sigma_i ) \Big) dx}  =
\begin{cases}
     0 & , \: \ i=N \\
     \int_{\partial \mathcal{D}}{U_D (\delta\sigma_i \cdot n) ds} & , \: \ i=D\\
\end{cases}
\label{eq:stateFluxes}
\end{equation}
Let $\sigma_{\Omega,N}^h$ and $\sigma_{\Omega,D}^h$ be the dual fluxes 
discretized using Raviart-Thomas Finite Element functions.
By combining the definition of energy-norm induced by the bilinear form 
(\ref{eq:a}) with the information in (\ref{eq:Optimality}) and 
(\ref{eq:stateFluxes}), we get the following upper bound for the energy-norm of 
the error in the state equations:
\begin{equation}
\enorm{u_{\Omega,i} - u_{\Omega,i}^h}^2 \leq \int_{\mathcal{D}}{\Big( 
k_\Omega^{-1} | \sigma_{\Omega,i}^h - k_\Omega \nabla u_{\Omega,i}^h|^2 + 
|\nabla \cdot \sigma_{\Omega,i}^h - u_{\Omega,i}^h|^2 \Big) dx}.
\label{eq:stateError}
\end{equation}

\subsubsection{Error estimates based on the complementary energy principle: the 
case of the adjoint equations}

As in the previous section, we present the formulation of the dual complementary 
energy associated with the discretization 
error of the adjoint problems (\ref{eq:adjointEIT}). 
In a similar fashion, we introduce the dual fluxes $\xi_{\Omega,i} \in \Hdiv$ 
for $i=N,D$ and we retrieve the following relationships
\begin{align}
\left\{
\begin{aligned}
& z = k_\Omega^{-1}\xi_{\Omega,i} + M(\theta^h)\nabla u_{\Omega,i}^h - \nabla 
r_{\Omega,i}^h \quad & \text{in} \ \mathcal{D} \\
& w = \nabla \cdot \xi_{\Omega,i} - \left( \nabla \cdot \theta^h \right) 
u_{\Omega,i}^h - r_{\Omega,i}^h\quad & \text{in} \ \mathcal{D} \\
& \xi_{\Omega,N} \cdot n = 0  \quad & \text{on} \ \partial\mathcal{D} 
\end{aligned}
\right.
\label{eq:OptimalityAdj}
\end{align}
The maximization problem for the dual complementary energy associated with the 
adjoint problems for $i=N,D$ reads as
\begin{equation}
\max_{\mytop{\xi_{\Omega,i} \in \Hdiv}{(\xi_{\Omega,N} 
\cdot n = 0)}}
-\frac{1}{2} \int_{\mathcal{D}}{\Big(k_\Omega^{-1} | \xi_{\Omega,i} + 
k_\Omega M(\theta^h)\nabla u_{\Omega,i}^h - k_\Omega \nabla 
r_{\Omega,i}^h|^2  + |\nabla \cdot \xi_{\Omega,i} - (\nabla \cdot 
\theta^h)u_{\Omega,i}^h - r_{\Omega,i}^h|^2 \Big) dx}.
\label{eq:maxSAdj}
\end{equation}
In order to compute the dual flux variables, we seek
$\xi_{\Omega,N},\ \xi_{\Omega,D} \in \Hdiv$ such that $\xi_{\Omega,N} \cdot n = 
0$ on $\partial \mathcal{D}$ satisfying $\forall \delta\xi_N,\ \delta\xi_D \in 
\Hdiv \ \text{such that} \ \delta\xi_N \cdot n=0 \ \text{on} \ \partial 
\mathcal{D}$
\begin{equation}
\int_{\mathcal{D}}{\Big(k_\Omega^{-1} \xi_{\Omega,i}\cdot\delta\xi_i + (\nabla 
\cdot \xi_{\Omega,i} )( \nabla \cdot \delta\xi_i) \Big) dx}  = 
\int_{\mathcal{D}}{\Big( (\nabla \cdot \theta^h)u_{\Omega,i}^h \nabla \cdot 
\delta\xi_i - M(\theta^h) \nabla u_{\Omega,i}^h \cdot \delta\xi_i \Big) dx }
\label{eq:adjointFluxes}
\end{equation}
and via their Raviart-Thomas Finite Element approximations $\xi_{\Omega,i}^h$'s, 
we derive the upper bound of the adjoint errors in the energy-norm
\begin{equation}
\enorm{r_{\Omega,i} -r_{\Omega,i}^h}^2 \leq 
\int_{\mathcal{D}}{\Big(k_\Omega^{-1} | \xi_{\Omega,i}^h + k_\Omega  
M(\theta^h) \nabla u_{\Omega,i}^h - k_\Omega \nabla r_{\Omega,i}^h|^2 + 
|\nabla \cdot \xi_{\Omega,i}^h - (\nabla \cdot \theta^h)u_i^h - 
r_{\Omega,i}^h|^2 \Big) dx}. 
\label{eq:adjError}
\end{equation}
Eventually, by plugging (\ref{eq:stateError}) and (\ref{eq:adjError}) into 
(\ref{eq:UBeit}), we are able to explicitly compute the upper bound 
$\overline{E}$ of the error in the shape gradient.

\section{Numerical results}
\label{ref:numerics}

We present some numerical results of the application of the Certified Descent 
Algorithm (CDA) to the problem of Electrical Impedance Tomography. 
We remark that the simulations presented in this paper are based on a mesh 
moving approach for the deformation of the domain (Algorithm 
\ref{scpt:shape-opt-adaptive} - step 7). Thus, the procedure does not allow for 
topological changes and the correct number of inclusions within the material has 
to be set at the beginning of the algorithm. 
In order to account for the nucleation of new inclusions or the merging of 
existing ones, a mixed approach based on topological and shape gradients may be 
followed as suggested in \cite{Hintermuller2008}.
For the rest of this section, we will consider several examples where the shape 
and the location of the inclusions evolve under the assumption that the number 
of subregions inside $\mathcal{D}$ is \emph{a priori} set and known. \\
It is well-known in the literature that the EIT problem is severely ill-posed.
Several approaches have been discussed in the literature, spanning shape 
optimization \cite{LaurainSturm2015, MR2211069, MR2407028}, 
topology optimization \cite{Hintermuller2008, MR2886190, MR2966180} and 
regularization methods \cite{MR3019471, MR2132313, holder2004electrical}.
Classical shape optimization methods are known to provide fairly poor 
results for the problem of Electrical Impedance Tomography, remaining 
trapped in local minima.
The Certified Descent Algorithm does not aim at solving this issue and 
results similar to the ones provided by the Boundary Variation Algorithm without 
certification are expected. Nevertheless, it is interesting to observe that 
an improved version of the classical Boundary Variation Algorithm still 
experiences issues handling the EIT problem.
We recall that the choice of the EIT problem is a proof of concept to establish some 
properties of the Certified Descent Algorithm on a non-trivial scalar case.
Moreover, the CDA acts as a counterexample confirming the limitations of 
gradient-based strategies when dealing with ill-posed problems as the Electrical 
Impedance Tomography.

Before running the shape optimization algorithm, we identify a set of consistent 
boundary conditions $(g,U_D)$ for the state problems (\ref{eq:stateEIT}).
First, we set a Neumann boundary condition $g$ on $\partial\mathcal{D}$ for the 
flux; in order to impose the Dirichlet condition on the potential, we 
iteratively solve the Neumann state problem by subsequently refining the mesh 
until a very fine error estimate in the energy-norm is achieved. 
The trace of the resulting solution $u_{\Omega,N}$ on $\partial\mathcal{D}$ is 
eventually picked as boundary datum $U_D$ for the Dirichlet state problem. 
\\
All the numerical results were obtained using FreeFem++ \cite{MR3043640}.
\begin{figure}[htb]
\centering
\includegraphics[width=0.4\columnwidth]{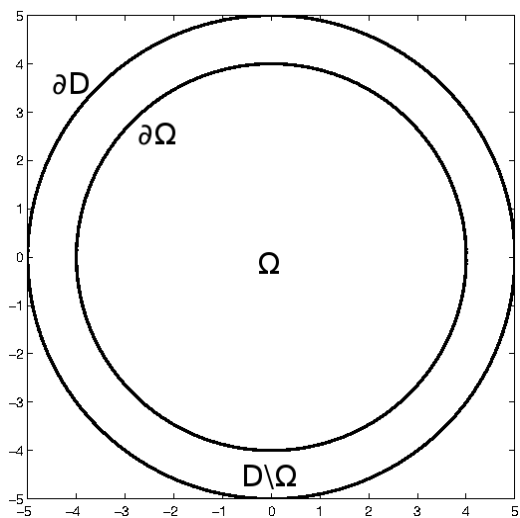}
\caption{Test case for the Electrical Impedance Tomography. Circular inclusion 
$\Omega$ inside the circular body $\mathcal{D}$.}
\label{fig:analyticalSol}
\end{figure}

\subsection{Numerical assessment of the goal-oriented estimator}
\label{ref:validation}

We consider the configuration in figure \ref{fig:analyticalSol}, where 
$\mathcal{D} \coloneqq \{ (x,y) \ | \ x^2+y^2 \leq \rho_E^2 \}$ 
and $\Omega \coloneqq \{ (x,y) \ | \ x^2+y^2 \leq \rho_I^2 \}$.
The values for the physical parameters are $\rho_E=5$, $\rho_I=4$, $k_E=1$ and 
$k_I=10$.
The boundary datum for the Neumann problem (\ref{eq:Neumann}) reads as $g=\cos(M 
\vartheta) \ , \ M=5$.
Using a polar coordinate system $(\rho,\vartheta)$, we can compute the following 
analytical solution: 
$$
u_{\Omega,N} = \begin{cases}
                C_0 J_M \left(-i \rho k_I^{-\frac{1}{2}}\right)\cos(M \vartheta) 
\ & , \ \rho \in [0,\rho_I]\\
                \left[ C_1 J_M \left(-i \rho k_E^{-\frac{1}{2}}\right) + C_2 Y_M 
\left(-i \rho k_E^{-\frac{1}{2}}\right) \right] \cos(M \vartheta) \ & , \ \rho 
\in (\rho_I,\rho_E]
               \end{cases}
$$
where $J_M(\cdot)$ and $Y_M(\cdot)$ respectively represent the first- and 
second-kind Bessel functions of order $M$. 
The constants $C_0,\ldots,C_2$ are detailed in table \ref{tab:constants}.
\begin{table}
\centering
\begin{tabular}[hbt]{| c || l | l |}
\hline
Constant & $\mathbb{R}\text{e}[C_i]$ & $\mathbb{I}\text{m}[C_i]$ \\
\hline & &
\\ [-1em] \hline
$C_0$ & $-6.3 \cdot 10^{-9}$ & $+40.39491005$ \\
\hline
$C_1$ & $+1.30145994$ & $+0.325482825$ \\
\hline
$C_2$ & $+1.5 \cdot 10^{-11}$ & $-1.301459935$ \\
\hline
\end{tabular}
\caption{Constants for the analytical solution.}
\label{tab:constants}
\end{table}

\begin{figure}[htb]
    \subfloat[Neumann problem.]
    {
    \includegraphics[width=0.45 \columnwidth]{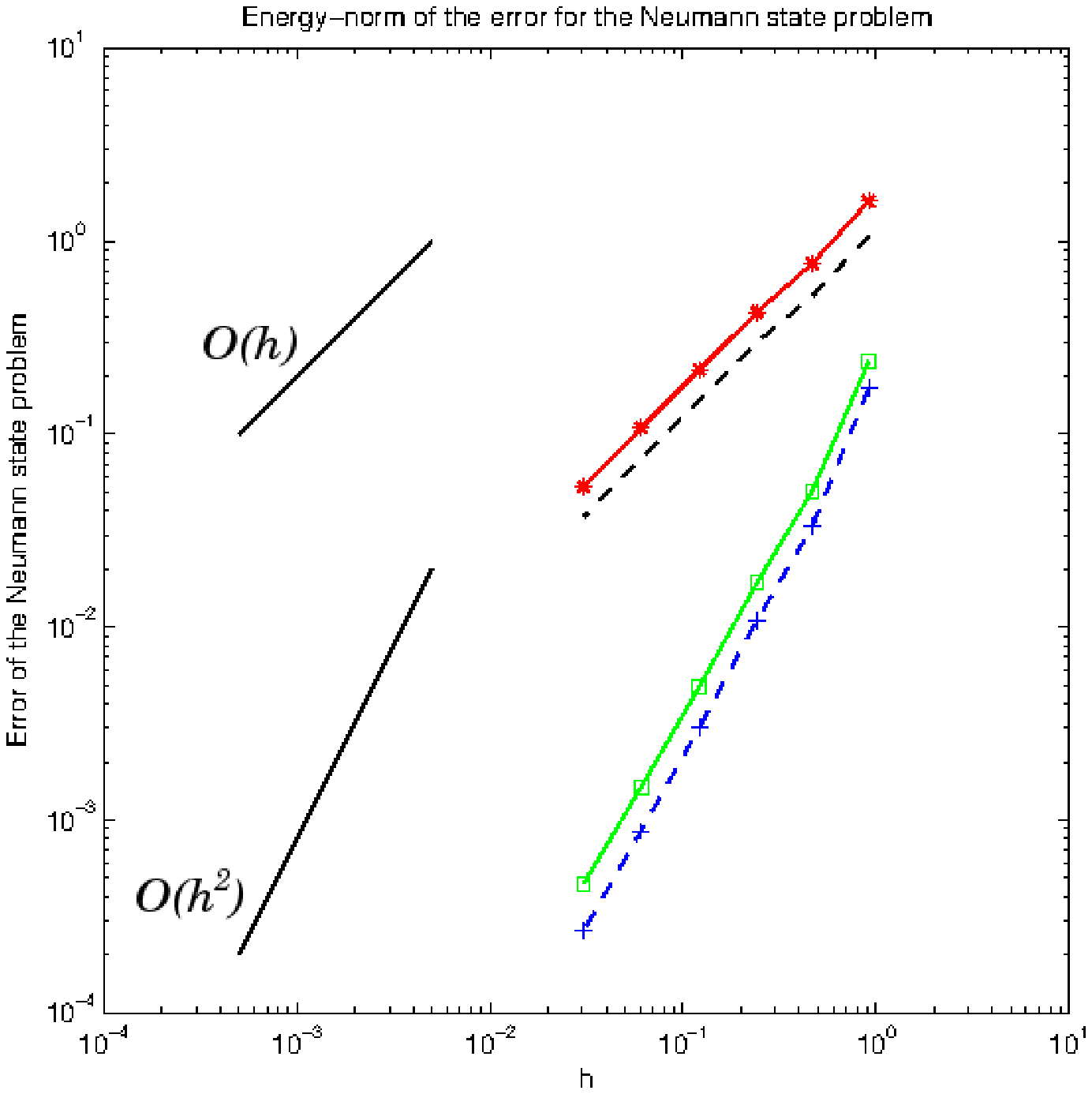}
    \vspace{10pt}
    \label{fig:stateN}
    }
    \hfil
    \subfloat[Dirichlet problem.]
    {
    \includegraphics[width=0.45 \columnwidth]{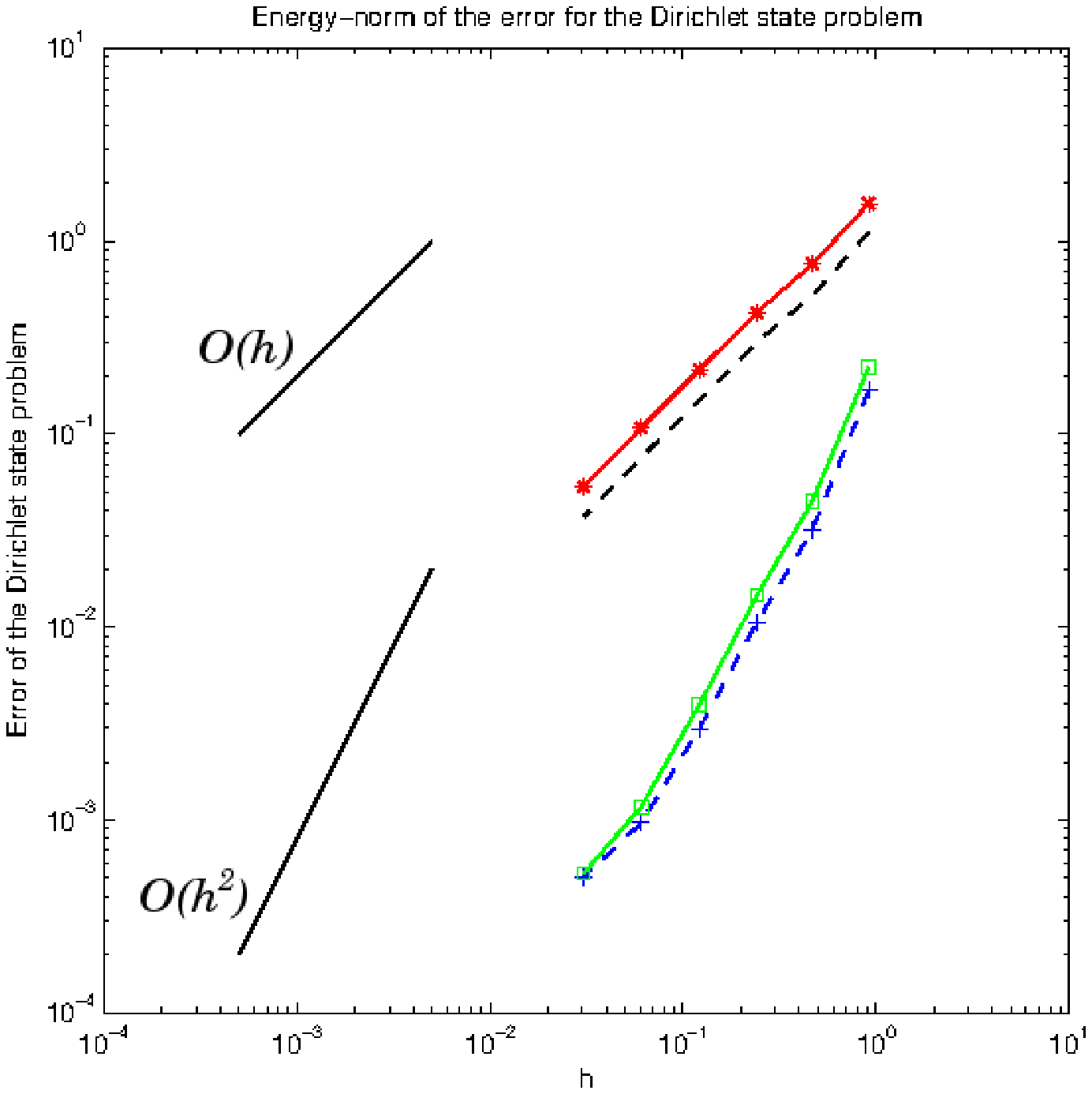}
    \label{fig:stateD}
    }
    \caption{Comparison of the convergence rates with respect to the mesh size 
$h$ for the (A) Neumann and (B) Dirichlet state equations. Analytical errors 
using $\mathbb{P}^1$ (black dash) and $\mathbb{P}^2$ (blue cross) Lagrangian 
Finite Element functions; error estimates based on the complementary energy 
principle using $\mathbb{P}^1$ (red star) and $\mathbb{P}^2$ (green square) 
Lagrangian Finite Element functions.}
    \label{fig:convergenceState}
\end{figure}

For the approximation of the state equations (\ref{eq:stateEIT}), we consider 
both $\mathbb{P}^1$ and $\mathbb{P}^2$ Lagrangian Finite Element functions. In 
figures \ref{fig:stateN} and \ref{fig:stateD}, we present a comparison between 
the analytical error due to the discretization and the corresponding estimates 
arising from the complementary energy principle (Section 
\ref{ref:stateEstimates}) under uniform 
mesh refinements.
We remark that using $\mathbb{P}^1$ Finite Elements, the estimated convergence 
rate is nearly $1$, whereas using $\mathbb{P}^2$ basis functions for the 
Finite Element space leads to a convergence rate slightly lower than $2$.\\
In order to construct the estimator for the error in the shape gradient, first 
we approximate equation (\ref{eq:appTheta}) using $\mathbb{P}^1 
\times \mathbb{P}^1$ Lagrangian Finite Element functions. For the 
discretization of the adjoint equations (\ref{eq:adjointEIT}), we use the same 
Finite Element space as for the state problems, whereas the dual fluxes in 
equations (\ref{eq:stateFluxes}) and (\ref{eq:adjointFluxes}) are approximated 
using Raviart-Thomas Finite Element functions. In particular, we choose the 
space of $RT_0$ (respectively $RT_1$) functions when the state and adjoint 
equations are solved using $\mathbb{P}^1$ (respectively $\mathbb{P}^2$) Finite 
Elements. 

\begin{figure}[hbt]
    \subfloat[$\mathbb{P}^1$ Finite Element.]
    {
    \includegraphics[width=0.45 \columnwidth]{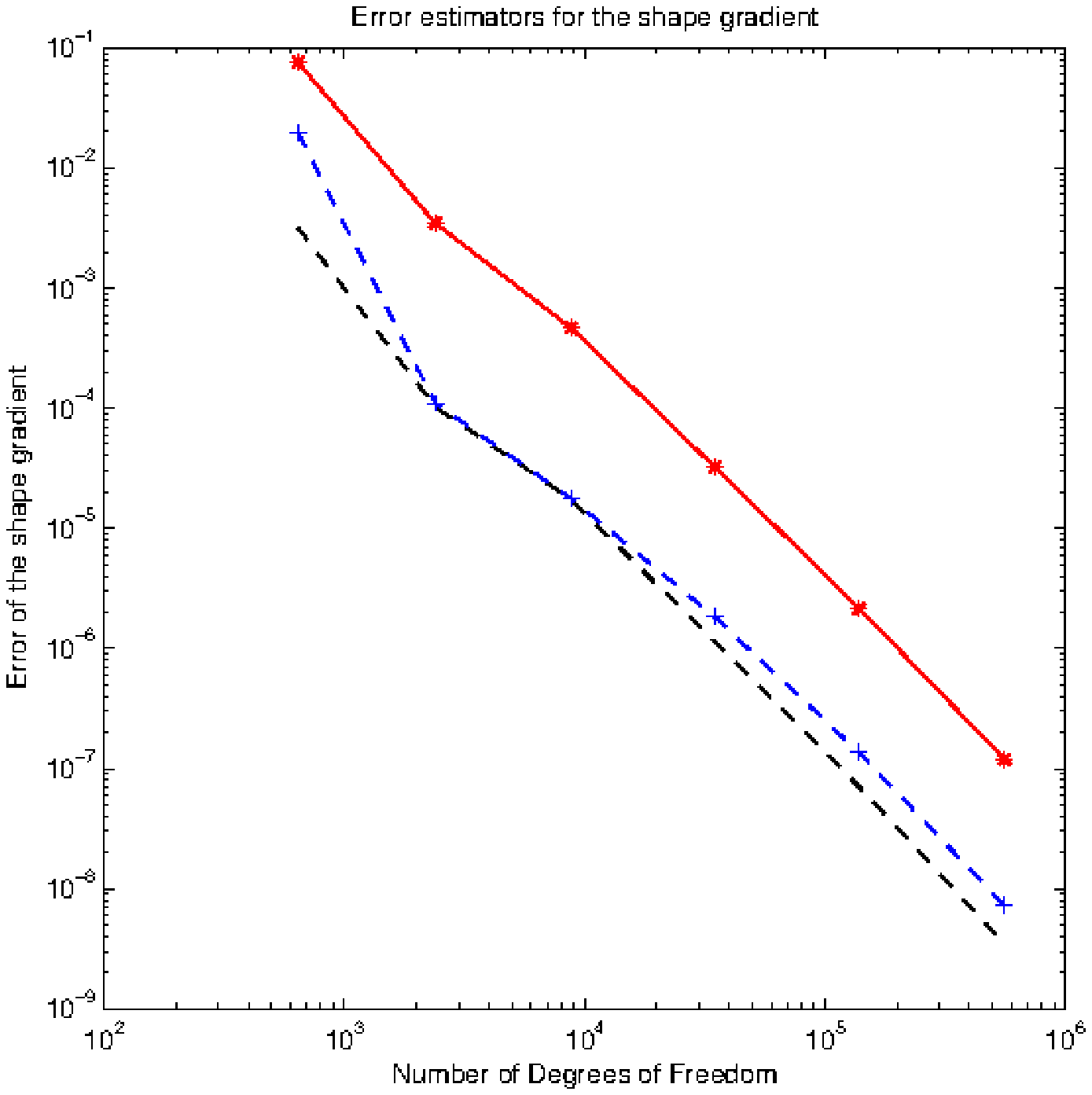}
    \vspace{10pt}
    \label{fig:QoIrateP1}
    }
    \hfil
    \subfloat[$\mathbb{P}^2$ Finite Element.]
    {
    \includegraphics[width=0.45 \columnwidth]{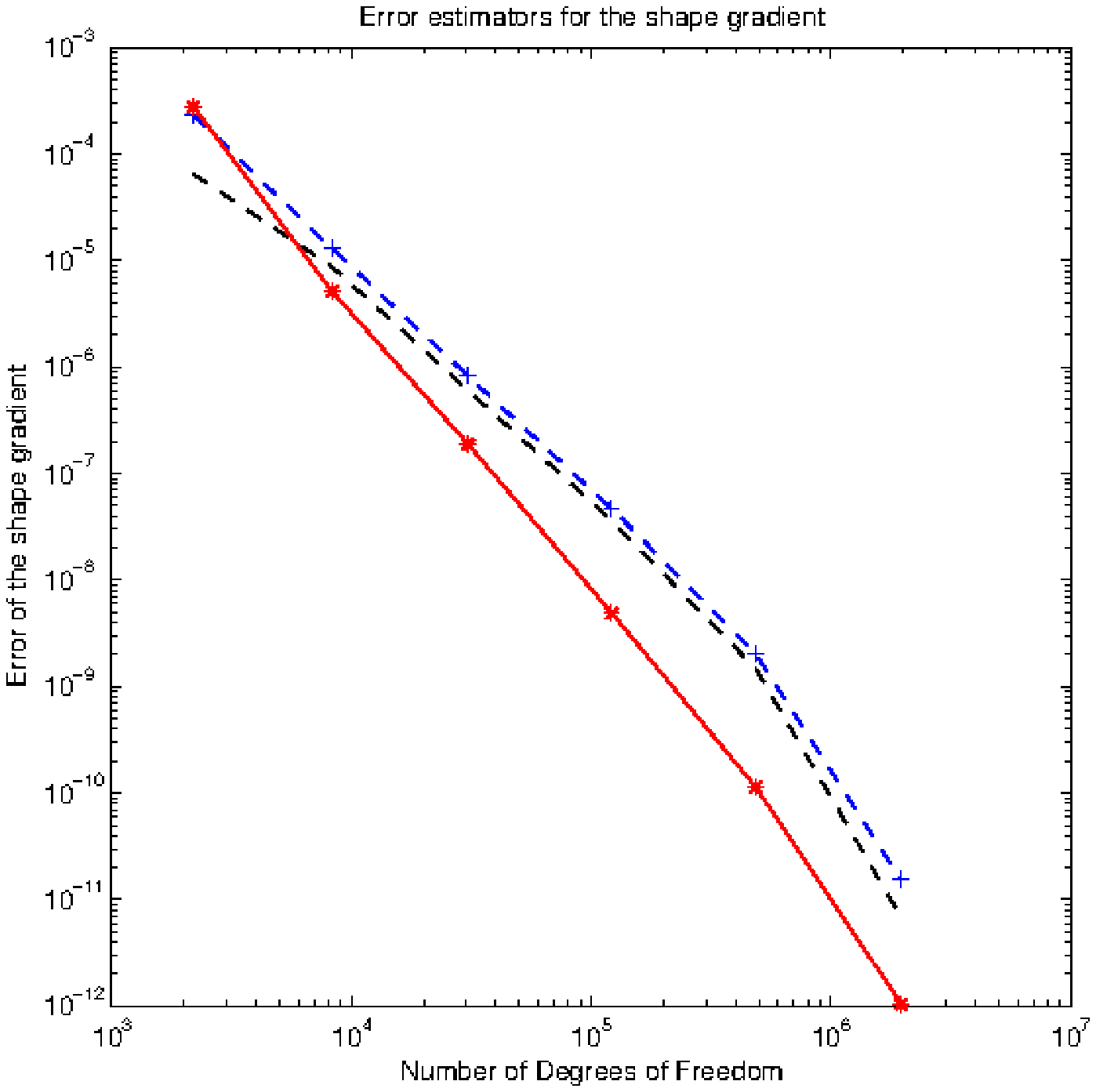}
    \label{fig:QoIrateP2}
    }
    \caption{Comparison of the convergence rates with respect to the number of 
Degrees of Freedom for the error in the Quantity of Interest 
    using (A) $\mathbb{P}^1$ and (B) $\mathbb{P}^2$ Lagrangian Finite Element 
functions.
    Analytical error in the shape gradient (black dash);
    error in the linearized Quantity of Interest (blue cross);
    error estimator for the Quantity of Interest using the complementary energy 
principle (red star).
    }
    \label{fig:QoIrate}
\end{figure}

In figure \ref{fig:QoIrate}, we present the convergence history of the 
discretization error in the shape gradient, the error in the Quantity of 
Interest arising from its linearization and the corresponding complementary 
energy estimates provided in (\ref{eq:UBeit}) using both $\mathbb{P}^1$ (Fig. 
\ref{fig:QoIrateP1}) and $\mathbb{P}^2$ (Fig. \ref{fig:QoIrateP2}) Finite 
Element functions.
In both figure \ref{fig:QoIrateP1} and figure \ref{fig:QoIrateP2}, we remark 
that the error in the linearized Quantity of Interest is very similar to the 
one in the shape gradient. 
This confirms that the linearization error introduced in (\ref{eq:errDJ}) is 
negligible with respect to the discretization error due to the Finite Element 
approximation and the estimator constructed from the linearized Quantity of 
Interest provides reliable information 
on the error in the shape gradient itself. \\
Figure \ref{fig:QoIrateP1} shows the evolution of the error in the Quantity of 
Interest with respect to the number of Degrees of Freedom of the problem under 
uniform mesh refinements. 
The error estimator shows an evolution which is analogous to the one of the 
analytical error in the shape gradient, thus we verify that an upper bound for 
the error in the Quantity of Interest is derived.
Nevertheless, when dealing with $\mathbb{P}^2$ Lagrangian Finite Element 
functions (Fig. \ref{fig:QoIrateP2}), the resulting error estimator for the 
Quantity of Interest underestimates the error in the shape gradient.
This phenomenon may be caused by the error due to the approximation of the 
geometry that has not been accounted for in this work. As a matter of fact, in 
\cite{COV:8787109} Morin \etal \ observe that increasing the accuracy of the PDE 
approximation is useless if the expected geometrical error is 
higher than the one due to the discretization of the state problem.
For this reason, in the following simulations, we stick to low-order Finite 
Element approximations ($\mathbb{P}^1-RT_0$) since using higher-order 
elements would prevent from getting a certified upper bound of the error in the 
shape gradient which is crucial for the application of the Certified 
Descent Algorithm.

\subsection{1-mesh and 2-mesh reconstruction strategies}

\begin{figure}[htbp]
    \centering
    \subfloat[$472$ elements.]
    {
    \includegraphics[width=0.3 \columnwidth]{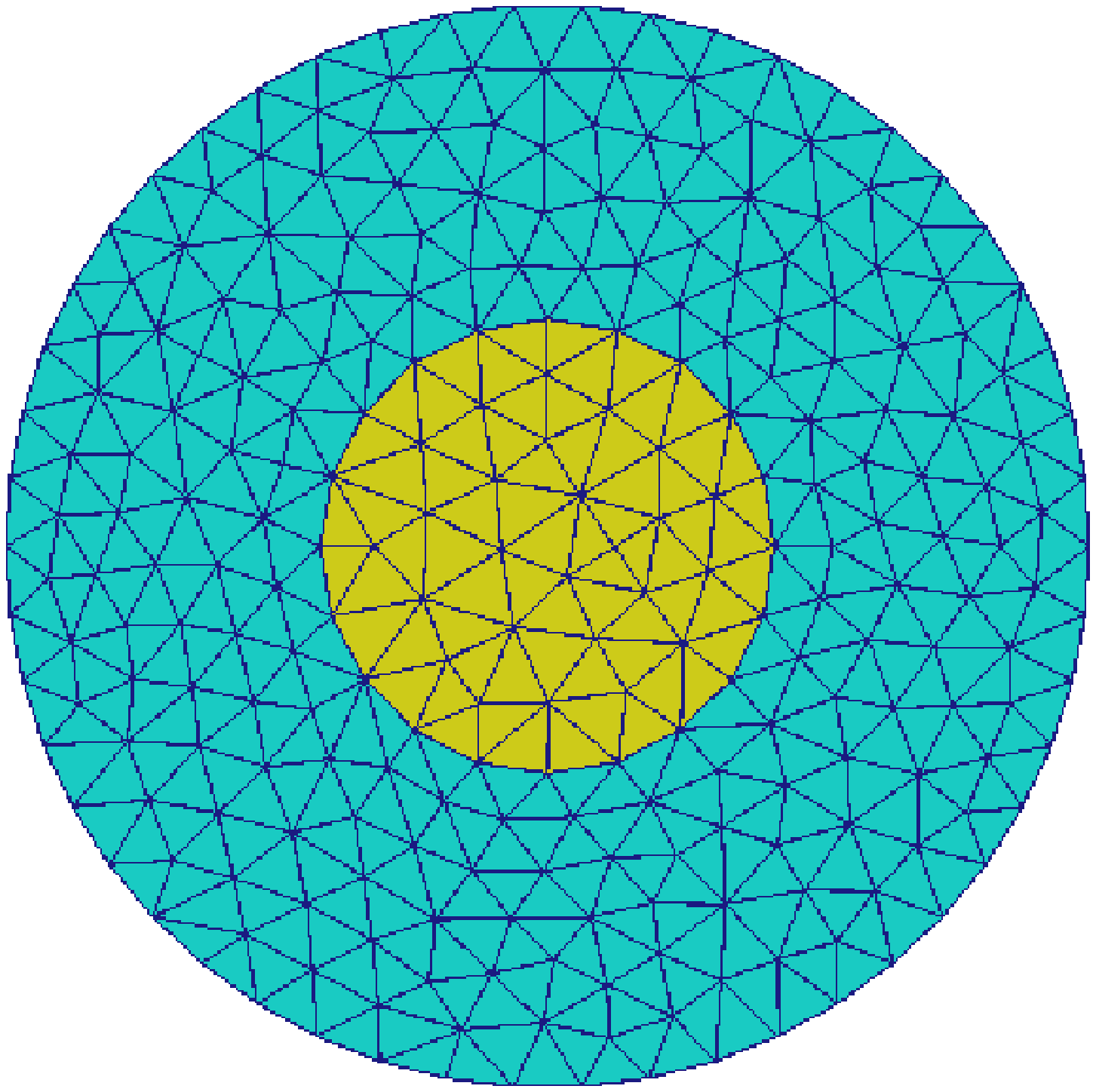}
    \label{fig:meshInCircle}
    }
    \hfil
    \subfloat[$178047$ elements.]
    {
    \includegraphics[width=0.3 \columnwidth]{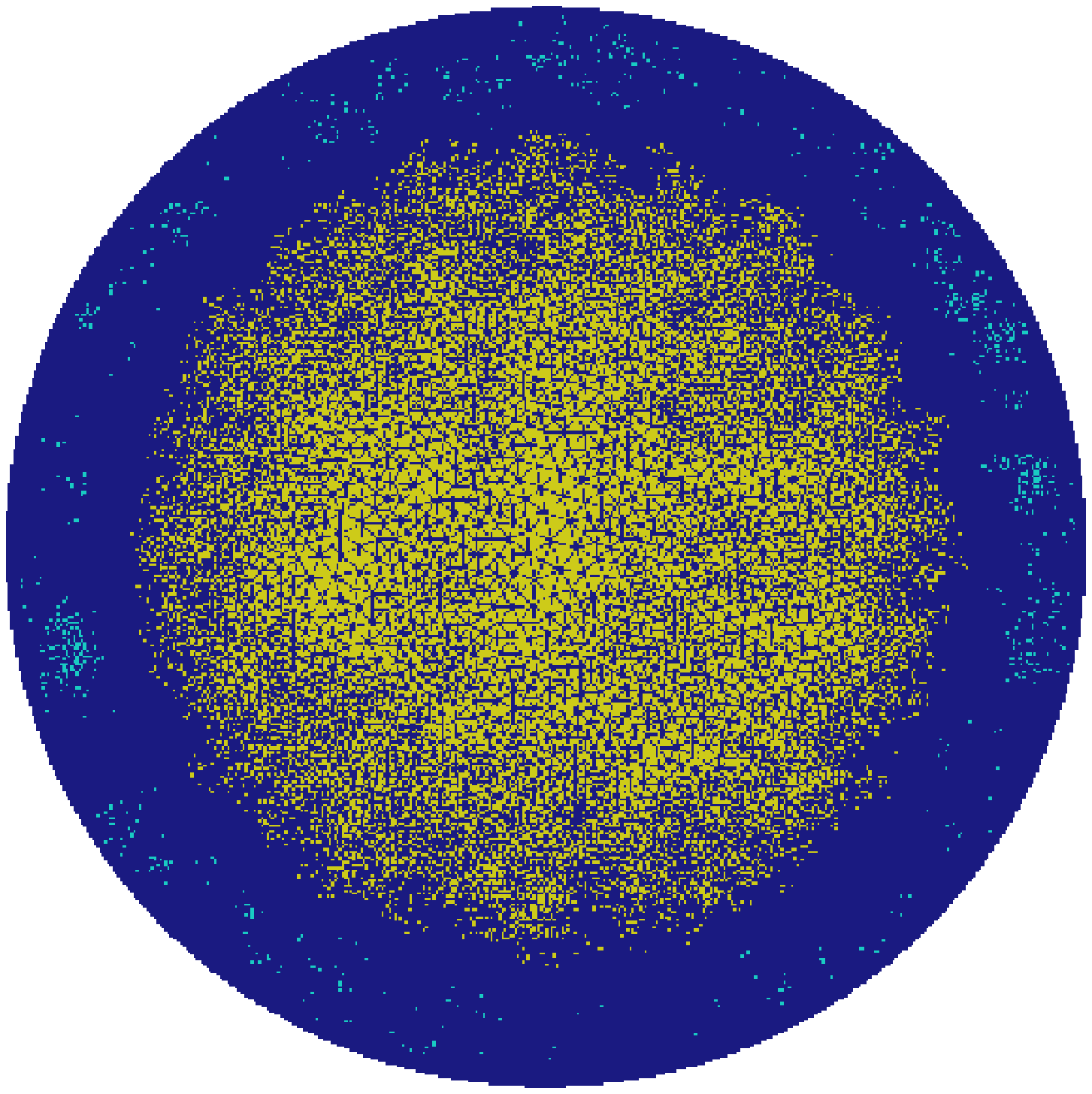}
    \label{fig:meshFinCircle1m}
    }
    \hfil
    \subfloat[$66833$ elements.]
    {
    \includegraphics[width=0.3 \columnwidth]{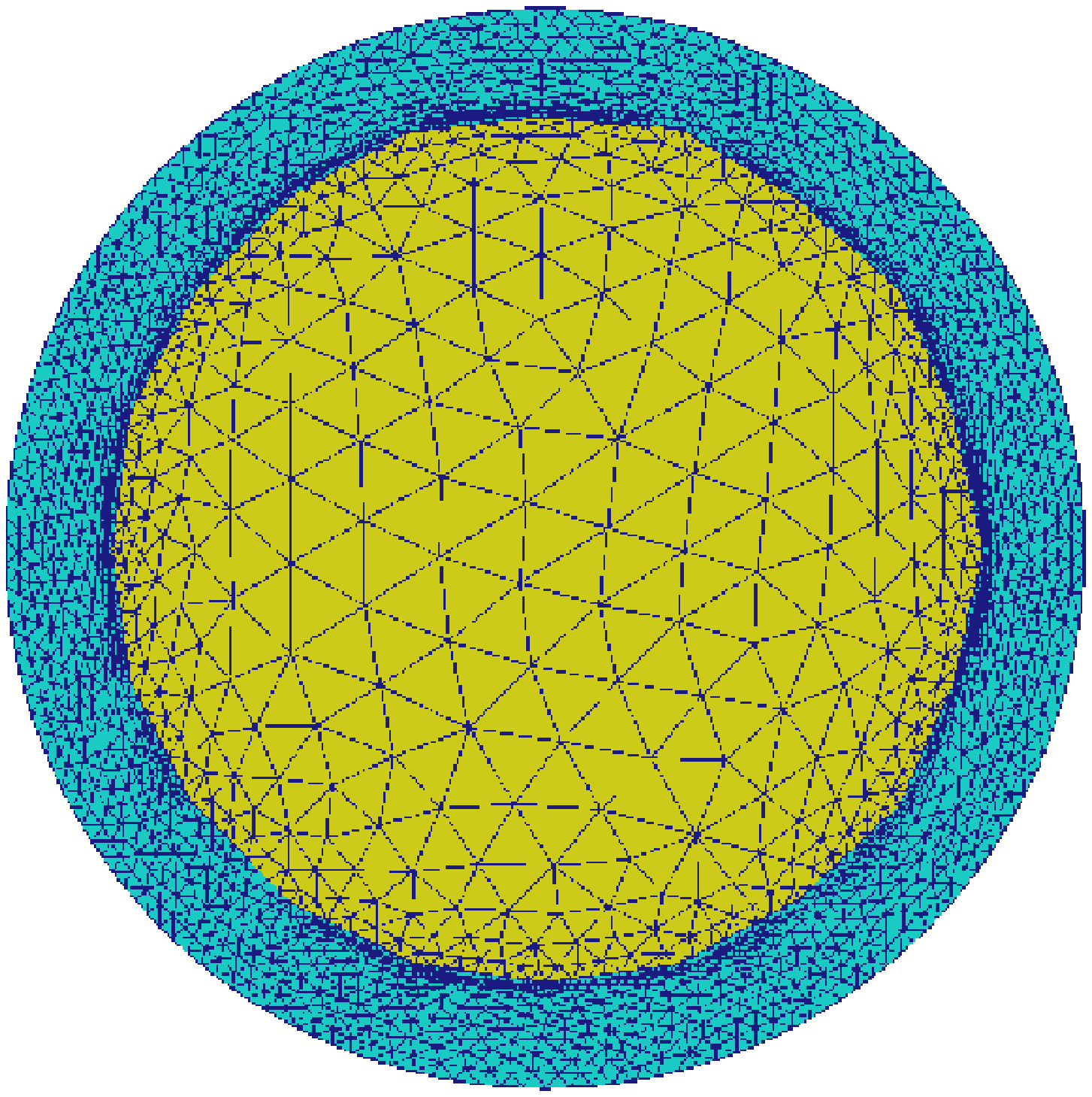}
    \label{fig:meshFinCircle2m}
    }
    
    \subfloat[Reconstructed interface.]
    {
    \includegraphics[width=0.3 \columnwidth]{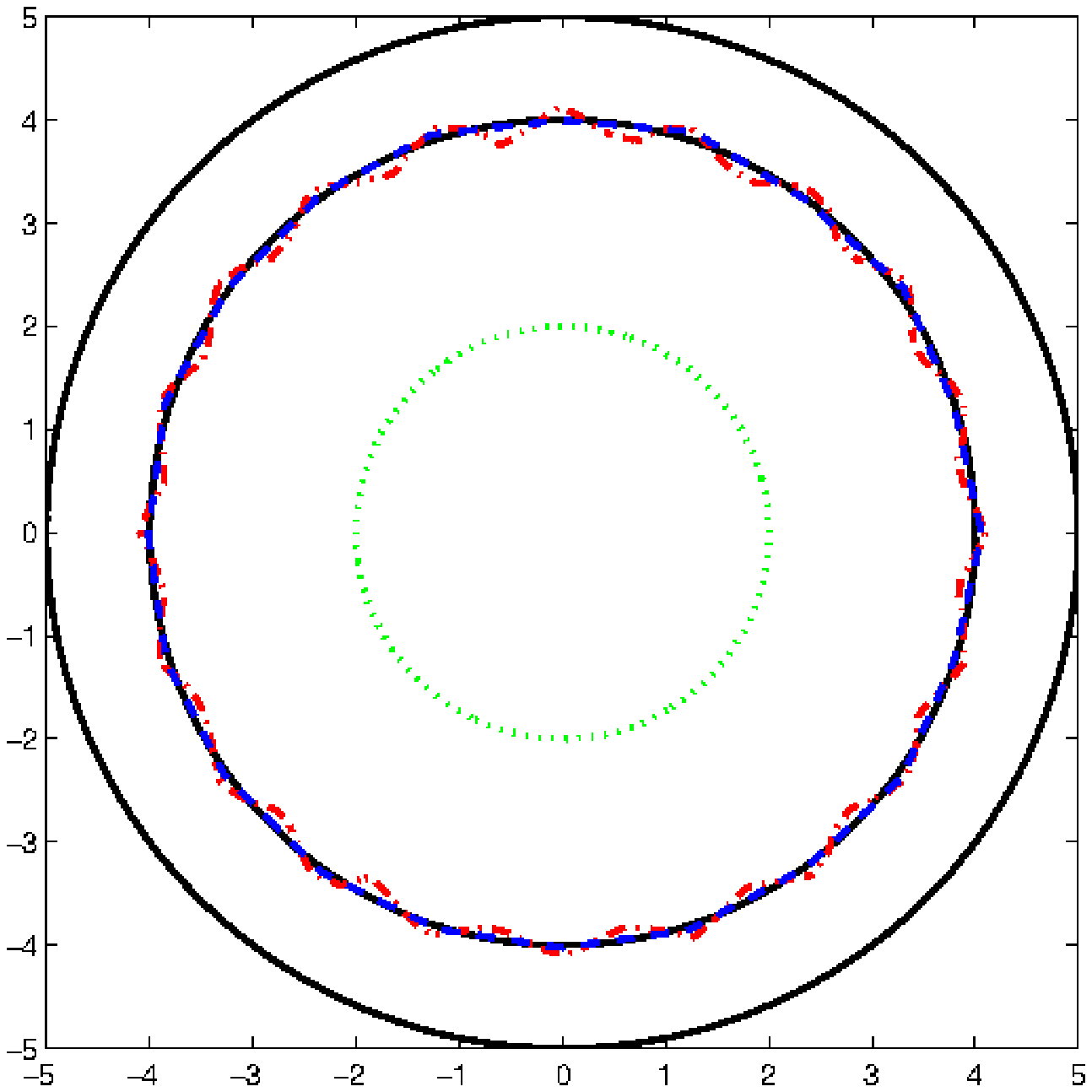}
    \label{fig:interfaceCircle}
    }
    \hfil
    \subfloat[1-mesh interface.]
    {
    \includegraphics[width=0.3 \columnwidth]{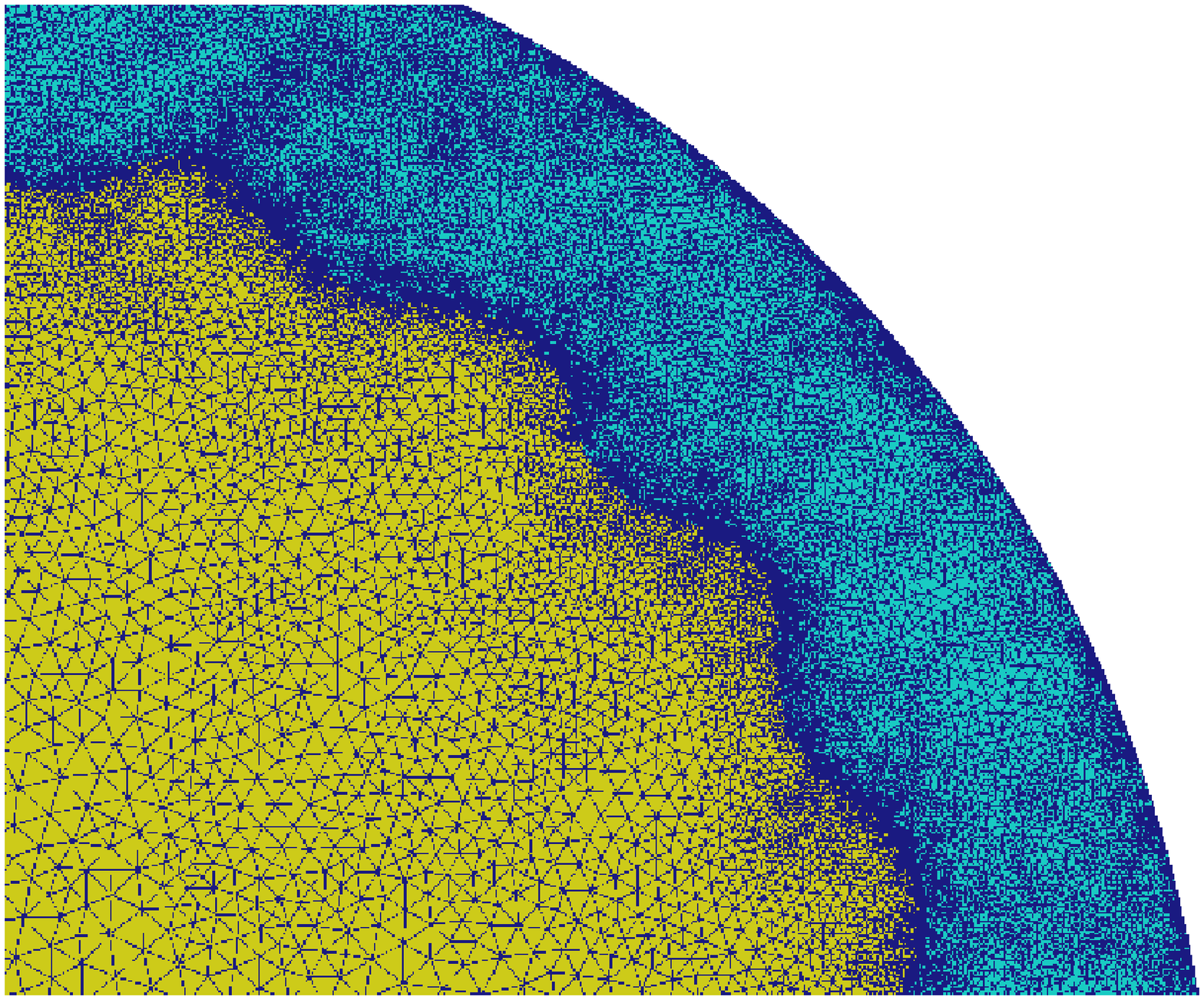}
    \label{fig:zoomCircle1m}
    }
    \hfil
    \subfloat[2-mesh interface.]
    {
    \includegraphics[width=0.3 \columnwidth]{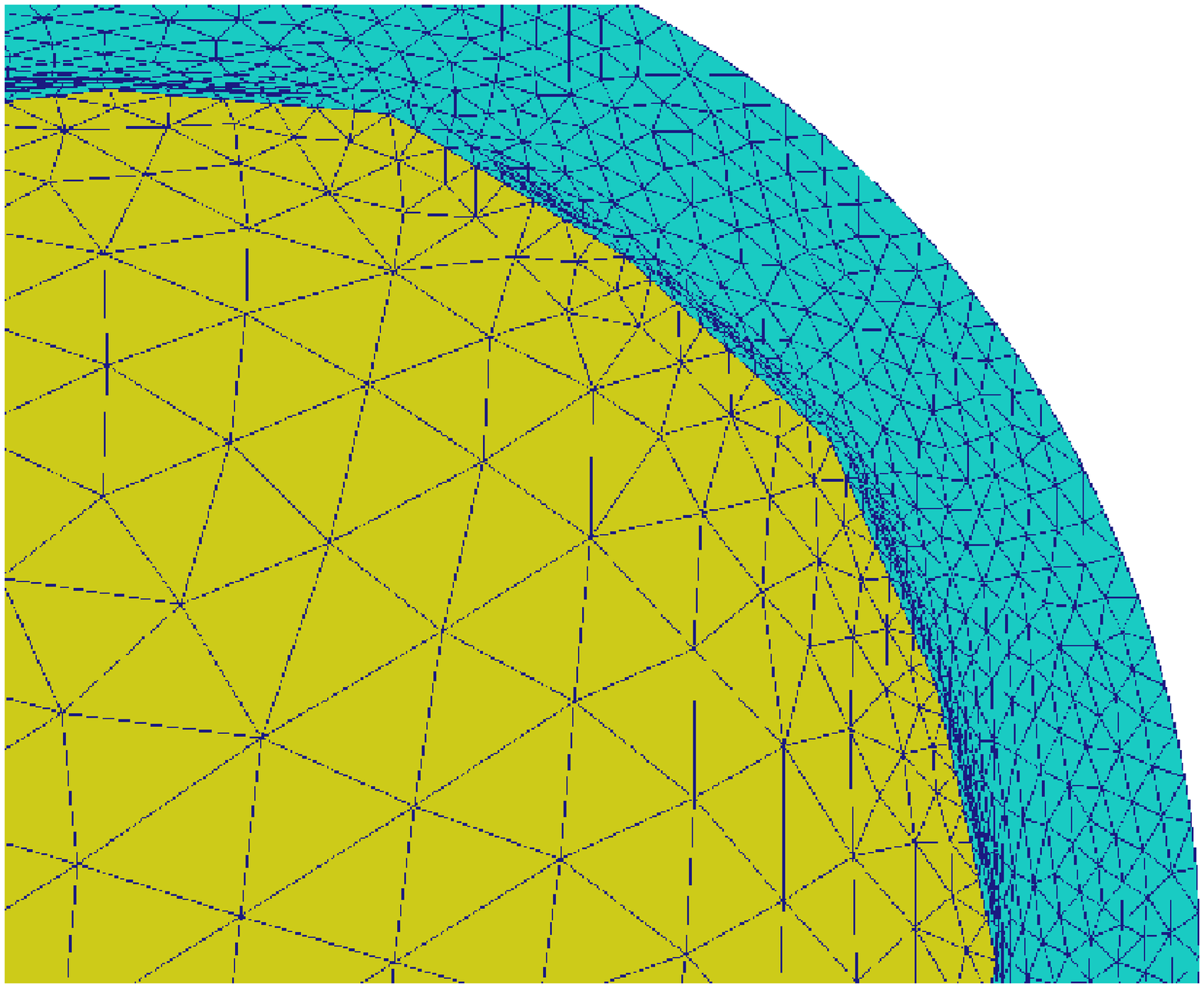}
    \label{fig:zoomCircle2m}
    }
    \caption{Comparison of 1-mesh and 2-mesh reconstruction strategies.
    (A) Initial mesh. (B) Final mesh using 1-mesh strategy. (C) Final mesh using 
2-mesh strategy. (D) Initial configuration (dotted green), target inclusion 
(solid black) and reconstructed interfaces using 1-mesh (dot-dashed red) and 
2-mesh (dashed blue) strategies. (E) Zoom of the reconstructed interface 
using 1-mesh strategy. (F) Zoom of the reconstructed interface using 2-mesh 
strategy.
    }    
    \label{fig:circleEx}    
\end{figure}

\begin{figure}[hbtp]
    \centering
    \subfloat[Objective functional.]
    {
    \includegraphics[width=0.48 \columnwidth]{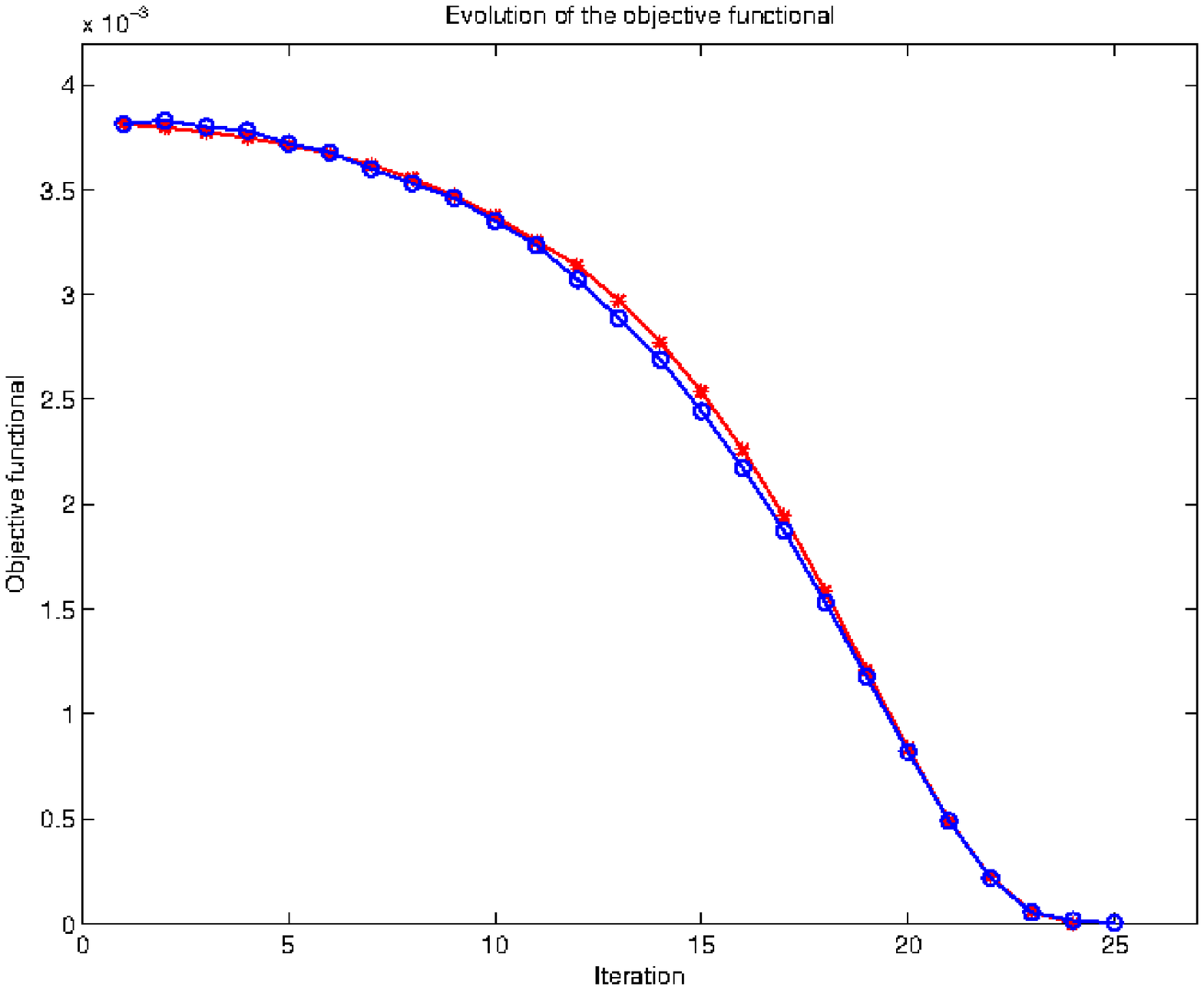}
    \label{fig:objCircle}
    }
    \hfil
    \subfloat[Number of Degrees of Freedom.]
    {
    \includegraphics[width=0.48 \columnwidth]{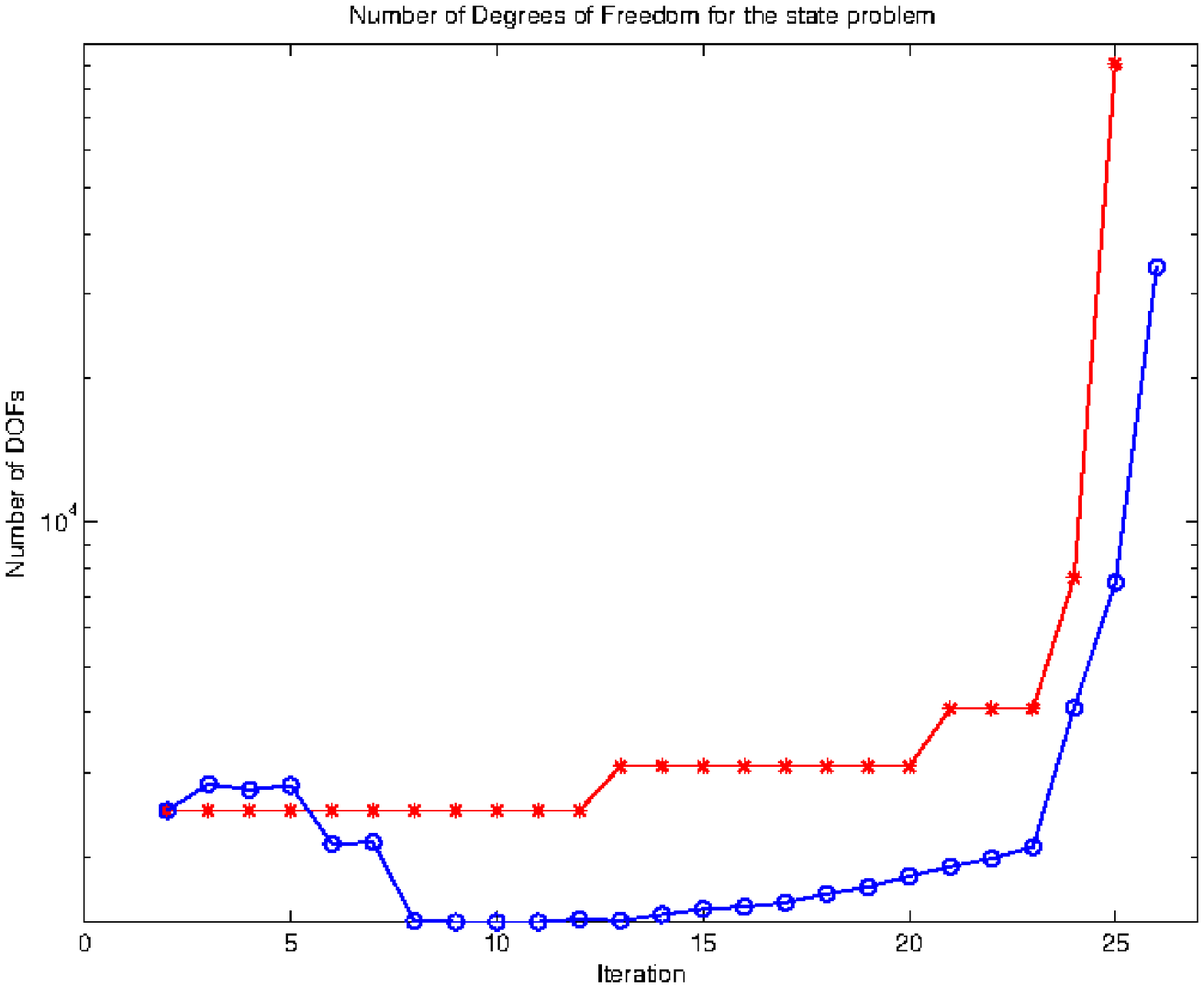}
    \label{fig:dofsCircle}
    }
    \caption{Comparison of 1-mesh (red star) and 2-mesh (blue circle) 
reconstruction strategies. (A) Evolution of the objective functional. (B) 
Number of Degrees of Freedom.
    }
    \label{fig:circleObjDofs}    
\end{figure}

We may now apply the CDA to identify the unknown inclusion $\Omega$ inside the 
circular domain $\mathcal{D}$ using one boundary measurement.
The initial inclusion is a circle of radius $\rho_{ini}=2$ and the associated 
computational mesh counting $472$ triangles is displayed in figure 
\ref{fig:meshInCircle}.\\
It is well-known in the literature (cf. \cite{smo-AP}) that using the same 
computational domain for both solving the state problem and computing a descent 
direction may lead to poor optimized shapes. 
In figures \ref{fig:meshFinCircle1m} and \ref{fig:meshFinCircle2m} we present 
the computational domains obtained using respectively a 1-mesh strategy to 
compute both the solutions $u_{\Omega,i}^h$'s and the descent direction 
$\theta^h$ and a 2-mesh approach in which the state equations are solved on a 
fine mesh whereas $\theta^h$ is computed on a coarser domain.
A comparison of the reconstructed interfaces after $24$ and $25$ iterations is 
reported in figure \ref{fig:interfaceCircle} and as expected we observe that 
the 1-mesh algorithm provides a poor approximation of the inclusion whereas the 
2-mesh strategy is able to precisely 
retrieve the boundary along which the conductivity $k_\Omega$ is discontinuous.
Figures \ref{fig:zoomCircle1m} and \ref{fig:zoomCircle2m} confirm what was 
already observed by zooming on the local behavior of the interfaces and 
highlighting the oscillatory nature of the 1-mesh reconstruction.
In figure \ref{fig:objCircle}, we report the evolution of the objective 
functional with respect to the number of iterations using the two discussed 
approaches. It is straightforward to observe that the CDA identifies a genuine 
descent direction at each iteration, generating a sequence of minimizing shapes 
such that the objective functional $J(\Omega)$ is monotonically decreasing. 
Moreover, the error estimate in the shape gradient is also 
used to construct a guaranteed stopping criterion for the overall optimization 
strategy which automatically ends when $| \langle d_h J(\Omega),\theta^h \rangle 
| + \overline{E} < \texttt{tol}$ for an \emph{a priori} set tolerance.
Even though both versions of the algorithm generate shapes for which the  
functional is monotonically decreasing, only the 2-mesh strategy is able to 
precisely identify the target inclusion. For this reason, in the following 
sections we will focus only on the 2-mesh approach. 

Besides the theoretical improvement of the Boundary Variation Algorithm 
discussed so far, an advantage of the CDA lies in the 
possibility of using relatively coarse meshes to identify certified descent 
directions. 
In figure \ref{fig:dofsCircle} we observe that the number of Degrees of Freedom 
remains small until the reconstructed interface approaches the real inclusion, 
that is coarse meshes prove to be reliable during the initial iterations of the 
algorithm.
Thus, another important feature of the Certified Descent Algorithm is the 
ability of certifying the reliability of coarse meshes for the identification of 
genuine descent directions for a shape functional, reducing the overall 
computational effort of the algorithm coupled with the \emph{a posteriori} 
estimators during the initial phase of computation.

\subsection{A more involved test case}
  
In the previous section, we applied the Certified Descent Algorithm to a simple  
test case and we were able to retrieve a precise description of the interface 
$\partial\Omega$. 
This was mainly due to the fact that the inclusion $\Omega$ was located near the 
external boundary $\partial \mathcal{D}$ where the measurements are performed. 
In this section, we consider a more involved test case: on the one hand, we 
break the symmetry of the problem by considering the initial and target 
configurations in figure \ref{fig:ellipse110}; on the other hand, we highlight 
the difficulties of precisely identifying the boundary of an inclusion when its 
position is far away from $\partial\mathcal{D}$.

As observed in the previous section, the evolution of the objective functional 
is monotonically decreasing (Fig. \ref{fig:objEllipse110} - red curve), meaning a genuine 
descent direction is identified at each iteration of the optimization 
procedure. 
The final value of the approximated objective functional is  
$\mathcal{O}(10^{-4})$, in agreement with the zero value in the analytical 
optimal configuration of the inclusion.\\
Moreover, the evolution of the number of Degrees of Freedom (Fig. 
\ref{fig:dofsEllipse110} - red curve) shows that coarse meshes prove to be reliable during 
initial iterations. The size of the problem remains small for several successive 
iterations but after few tens of iterations the CDA performs multiple mesh refinements 
in order to identify a genuine descent direction. This results in a 
high number of Degrees of Freedom which rapidly increases when approaching the 
configuration for which the criterion $| \langle d_h J(\Omega),\theta^h \rangle | + \overline{E}$ 
fulfills a given tolerance $\texttt{tol}=5 \cdot 10^{-8}$. 

In figure \ref{fig:ellipse110}, we observe that the part of the interface closest 
to $\partial\mathcal{D}$ is well identified using the Certified Descent 
Algorithm with one boundary measurement. Nevertheless when moving away from the 
external boundary, the precision of the reconstructed interface decreases and 
the algorithm is not able to precisely identify the whole inclusion.
The uncertainty of the reconstruction in the central region of $\mathcal{D}$ is 
mainly due to the ill-posedness of the inverse problem. 
As a matter of fact, state problems (\ref{eq:Neumann}) and (\ref{eq:Dirichlet}) 
are elliptic equations thus the effect of the boundary conditions becomes less 
and less important moving away from $\partial\mathcal{D}$.

\subsubsection*{The case of multiple boundary measurements}

\begin{figure}[htb]
    \subfloat[Reconstructed interface.]
    {
    \includegraphics[width=0.28 \columnwidth]{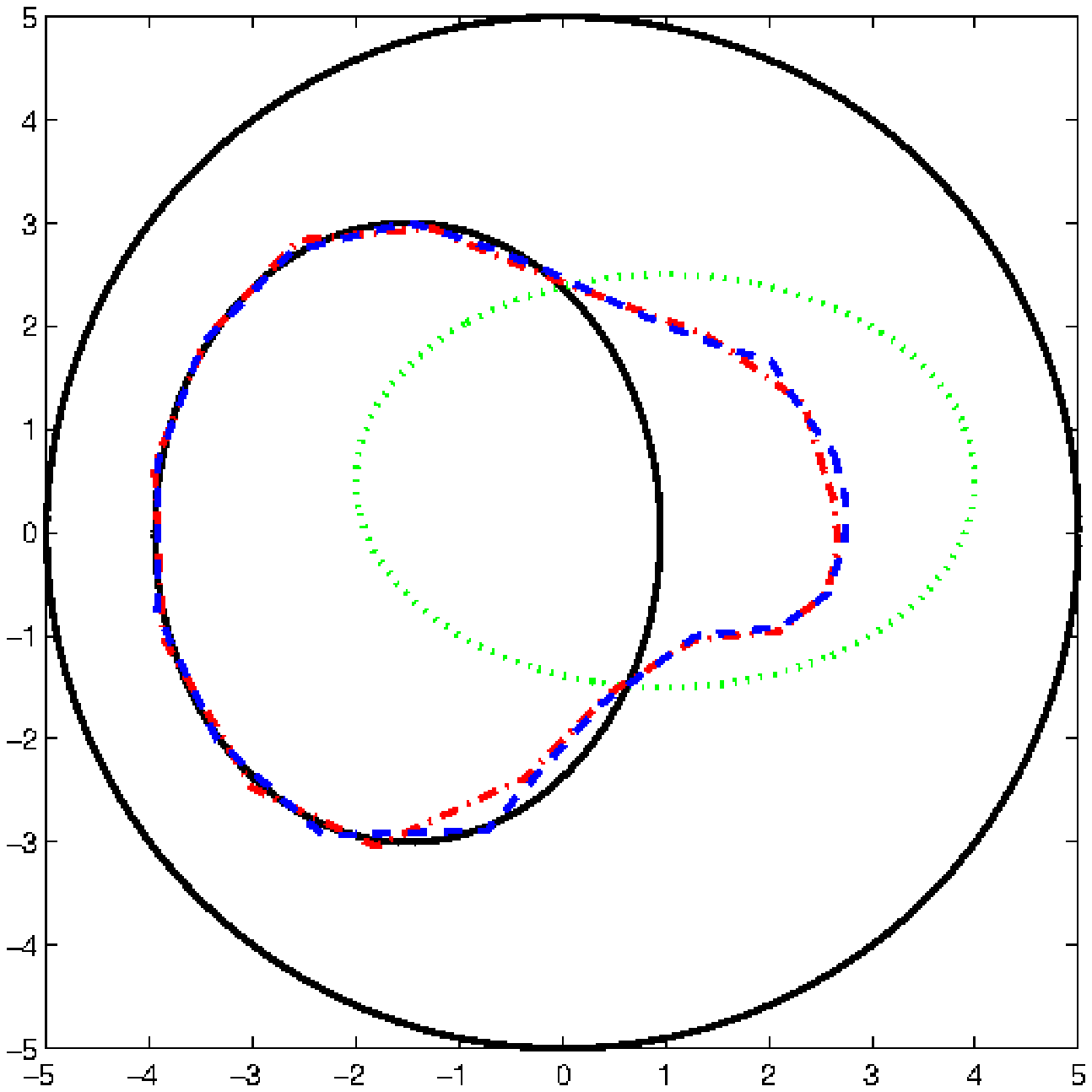}
    \vspace{10pt}
    \label{fig:ellipse110}
    }
    \hfil
    \subfloat[Objective functional.]
    {
    \includegraphics[width=0.33 \columnwidth]{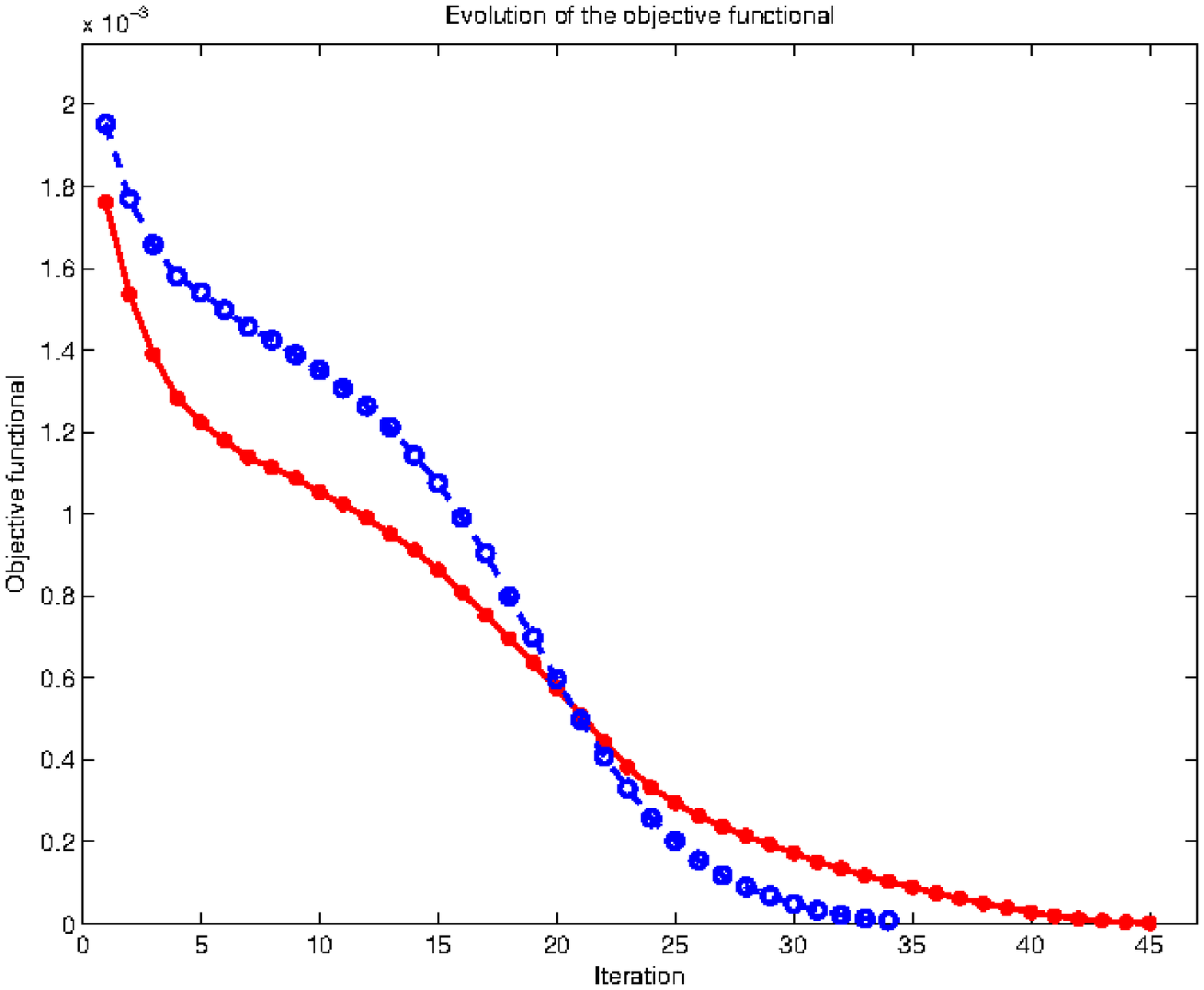}
    \label{fig:objEllipse110}
    }
    \hfil
    \subfloat[Degrees of Freedom.]
    {
    \includegraphics[width=0.33 \columnwidth]{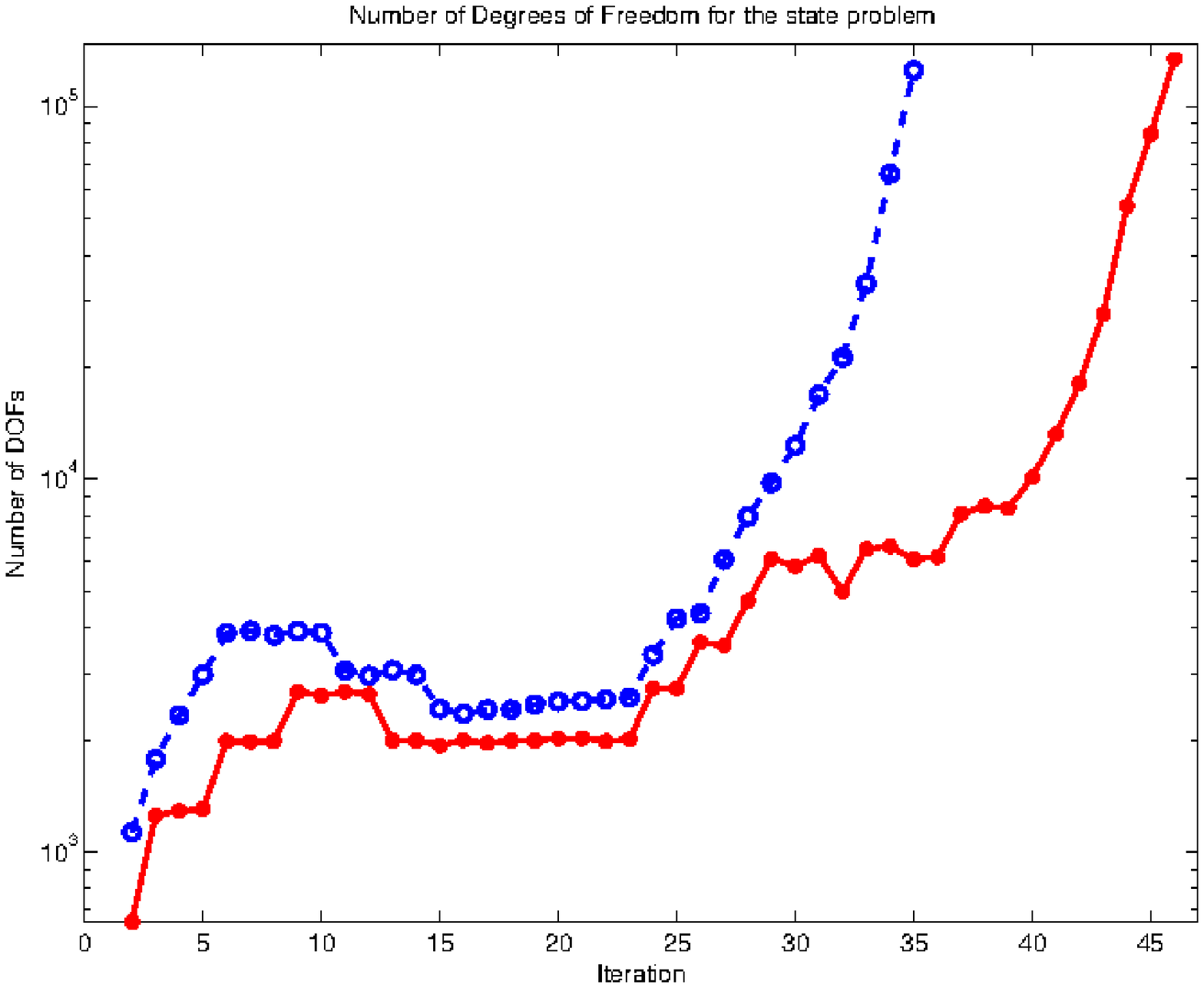}
    \label{fig:dofsEllipse110}
    }
    
    \caption{Certified Descent Algorithm using one measurement (red stars) and ten 
    measurements (blue circles). 
    (A) Initial configuration (dotted green), target inclusion (solid black) and 
reconstructed interface (dot-dashed red and dashed blue). (B) Evolution of the objective 
functional. (C) Number of Degrees of Freedom.}
    \label{fig:ellipse}
\end{figure}

A strategy to improve the quality of the reconstruction via the CDA relies 
on the use of several boundary measurements to retrieve a better approximation 
of the inclusion in a smaller number of iterations. 
The procedure to construct the set of boundary conditions to be used by the 
CDA is detailed in next section. Here, we present the outcome of the 
Certified Descent Algorithm using ten measurements for the test case 
featuring the circular domain $\mathcal{D} \coloneqq \{ (x,y) \ | \ x^2+y^2 \leq 25 \}$ 
with the inclusion represented by the solid black line in figure \ref{fig:ellipse110}.
We observe that using several boundary measurements the algorithm has access 
to more information to better identify the interface of the inclusion. 
First of all, the blue curve in figure \ref{fig:ellipse110} confirms the ability 
of the CDA to exactly identify the inclusion near the boundary $\partial\mathcal{D}$ 
and it highlights some minor improvements in the reconstruction of the interface 
with respect to the test case in red featuring one measurement 
(cf. the upper and lower parts of the ellipse).
Nevertheless, the result is still degraded when moving towards the center of the domain.
A possible workaround to this issue and to the low resolution of the 
reconstruction in the center of the computational domain is proposed by Ammari 
\etal \ in \cite{doi:10.1137/120863654}, where a hybrid imaging method arising 
from the coupling of electromagnetic tomography with acoustic waves is 
described. \\
The observed phenomenon is due to the well-known ill-posedness of the problem. 
Though considering several measurements improve the overall outcome of the algorithm, 
the final result is far from being satisfactory. Nevertheless, these limitations 
are related to the nature of the problem and we cannot expect gradient-based 
strategies to successfully overcome this issue.
These remarks are confirmed by figure \ref{fig:dofsEllipse110}. 
As a matter of fact, the number of Degrees of Freedom rapidly increases in both 
test cases, reaching $10^5$ and making the certification procedure unfeasible.
Besides the improvement in the reconstructed interface, 
the use of several measurements is responsible for reducing the number of iterations 
required by the CDA to identify the inclusion (Fig. \ref{fig:objEllipse110}). Moreover we remark that 
the tolerance that the quantity $| \langle d_h J(\Omega),\theta^h \rangle | + \overline{E}$ 
has to fulfill in this case drops to $\texttt{tol}=10^{-6}$, that is 
finer results are obtained in a smaller number of iterations and 
using lower precision in the case of multiple 
boundary measurements.

\subsubsection*{Sensitivity of the CDA to different initializations}

\begin{figure}[htb]
    \centering
    \subfloat[Initialization (I).]
    {
    \includegraphics[width=0.3 \columnwidth]{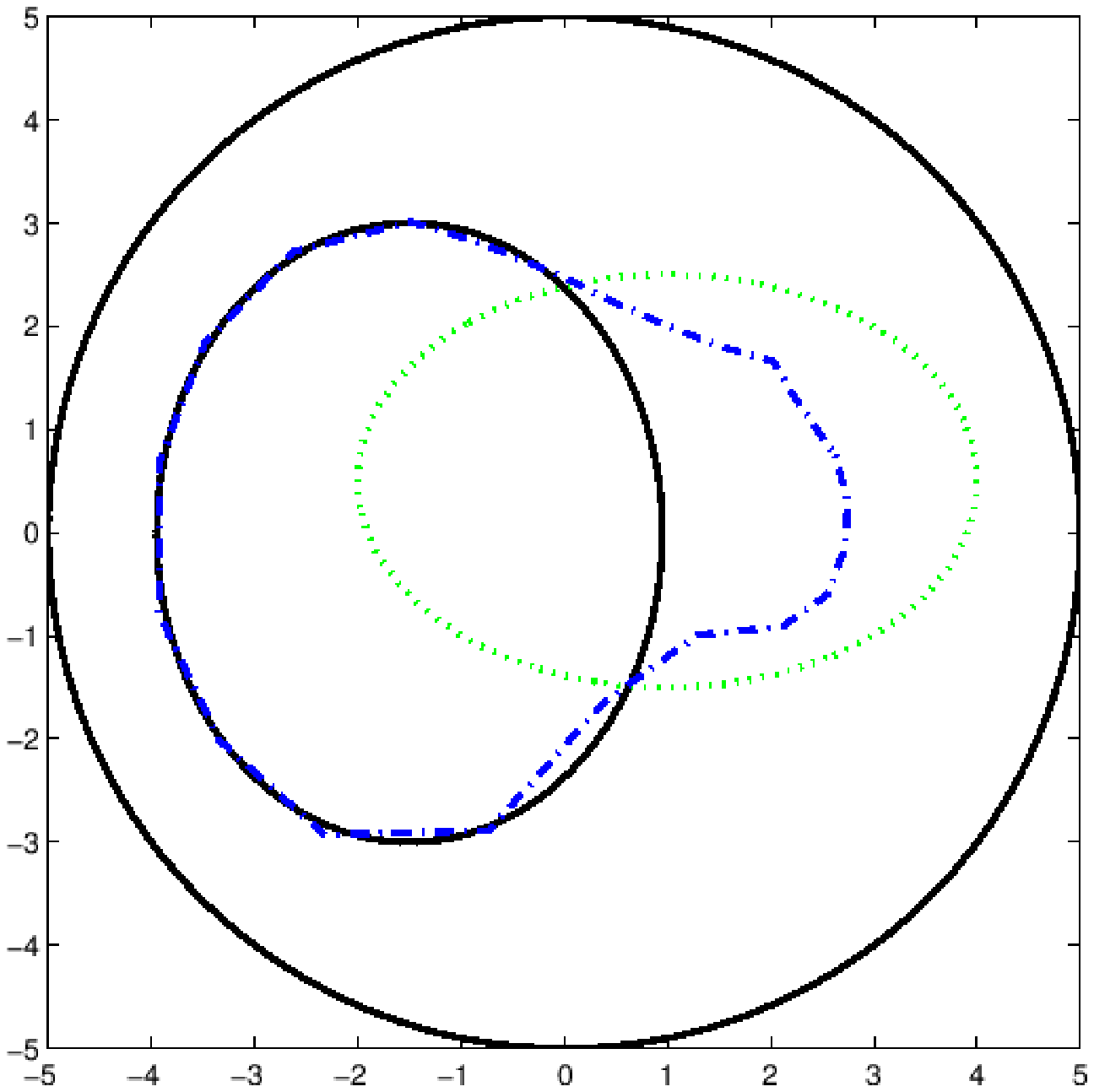}
    \label{fig:interfaceEll10}
    }
    \hfil
    \subfloat[Initialization (II).]
    {
    \includegraphics[width=0.3 \columnwidth]{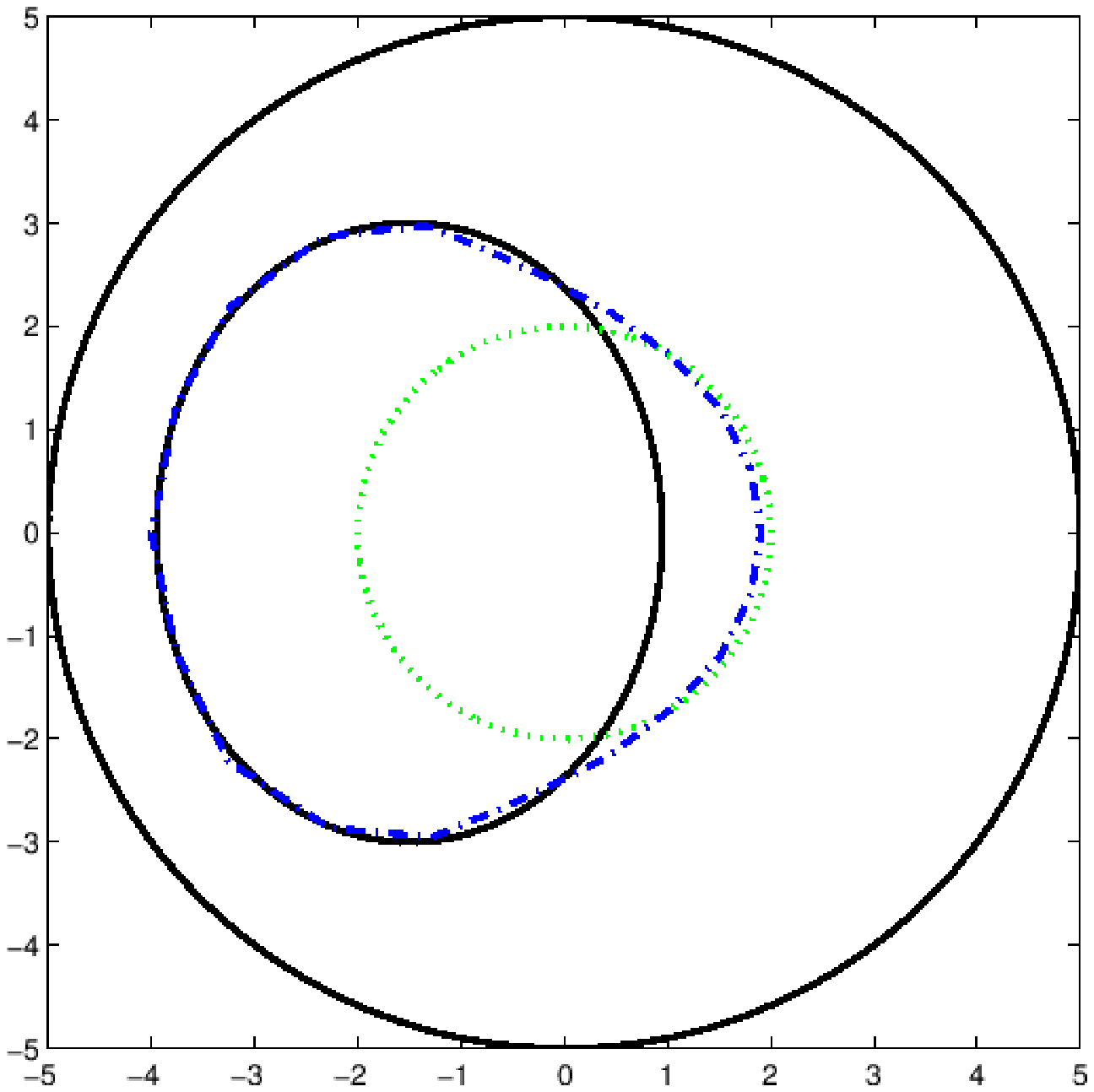}
    \label{fig:interfaceEllA}
    }
    \hfil
    \subfloat[Initialization (III).]
    {
    \includegraphics[width=0.3 \columnwidth]{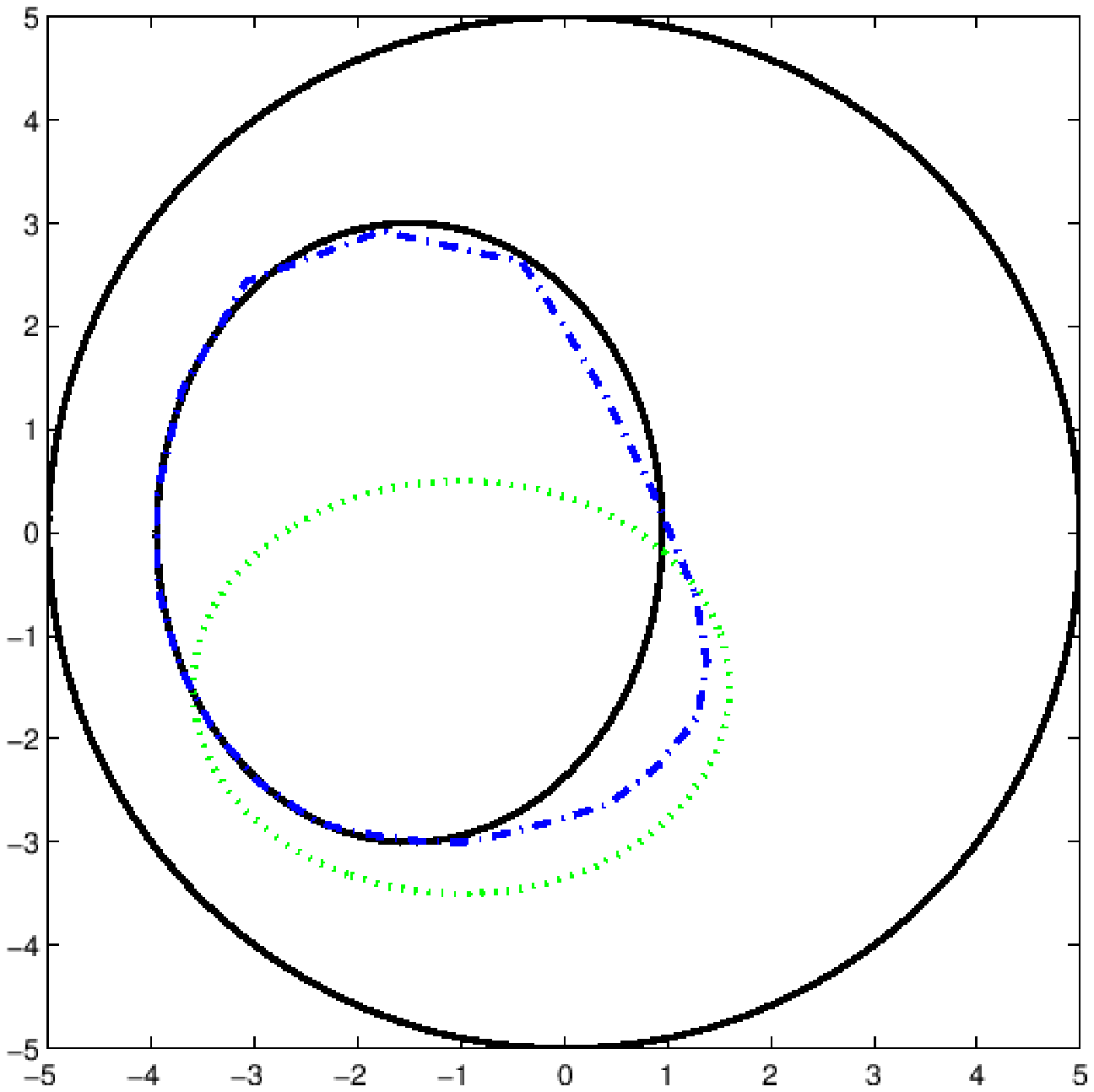}
    \label{fig:interfaceEllB}
    }

    \caption{Certified Descent Algorithm using multiple measurements 
    for different initial configurations.
    Initial configuration (dotted green), target inclusion (solid black) and reconstructed interface 
(dashed blue). 
    }
    \label{fig:ellipseInterfaceInitializations}    
\end{figure}

\begin{figure}[htb]
\centering
    \subfloat[Objective functional.]
    {
    \includegraphics[width=0.48\columnwidth]{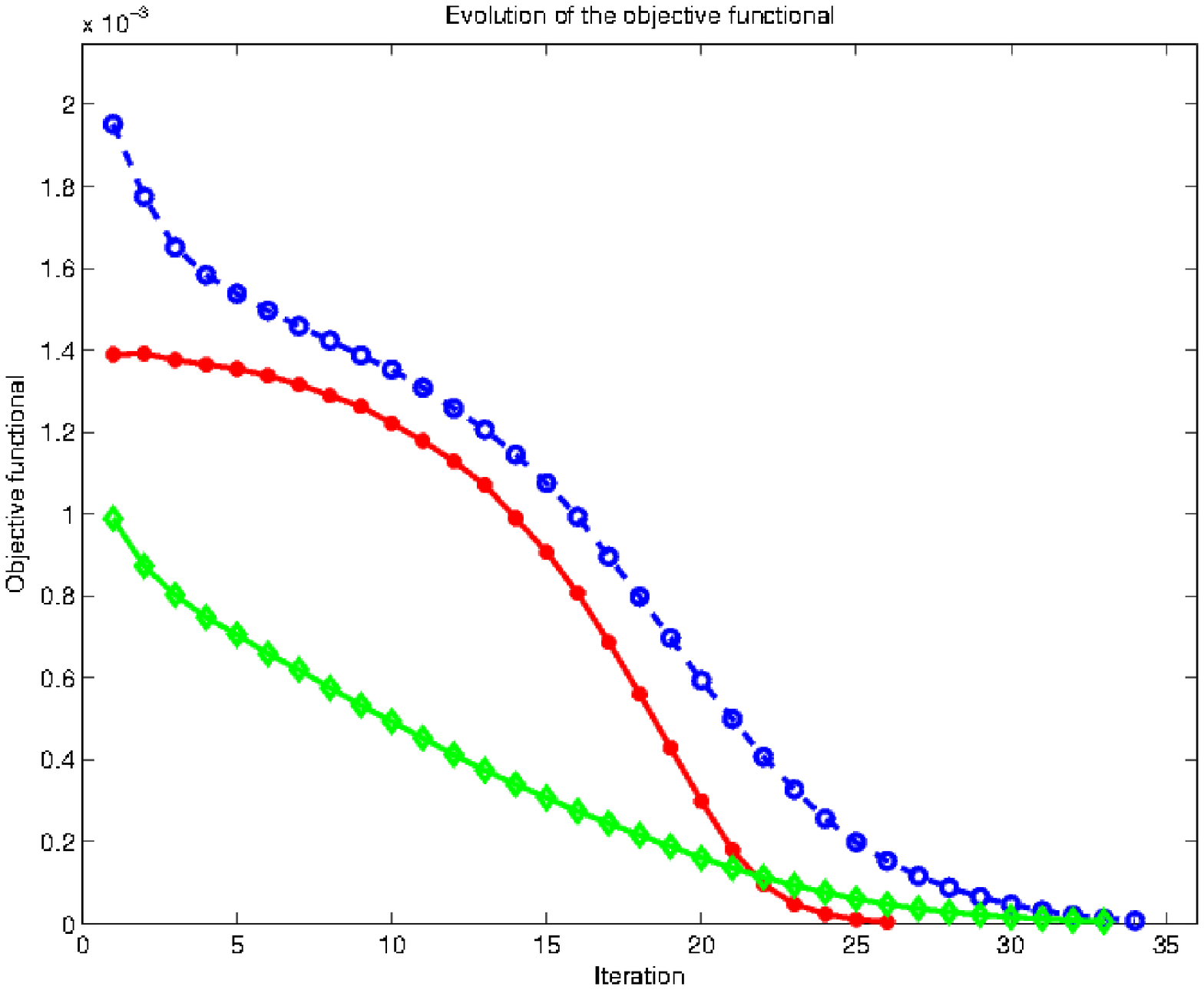}
    \label{fig:objEllInit}
    }
    \hfil
    \subfloat[Number of Degrees of Freedom.]
    {
    \includegraphics[width=0.48\columnwidth]{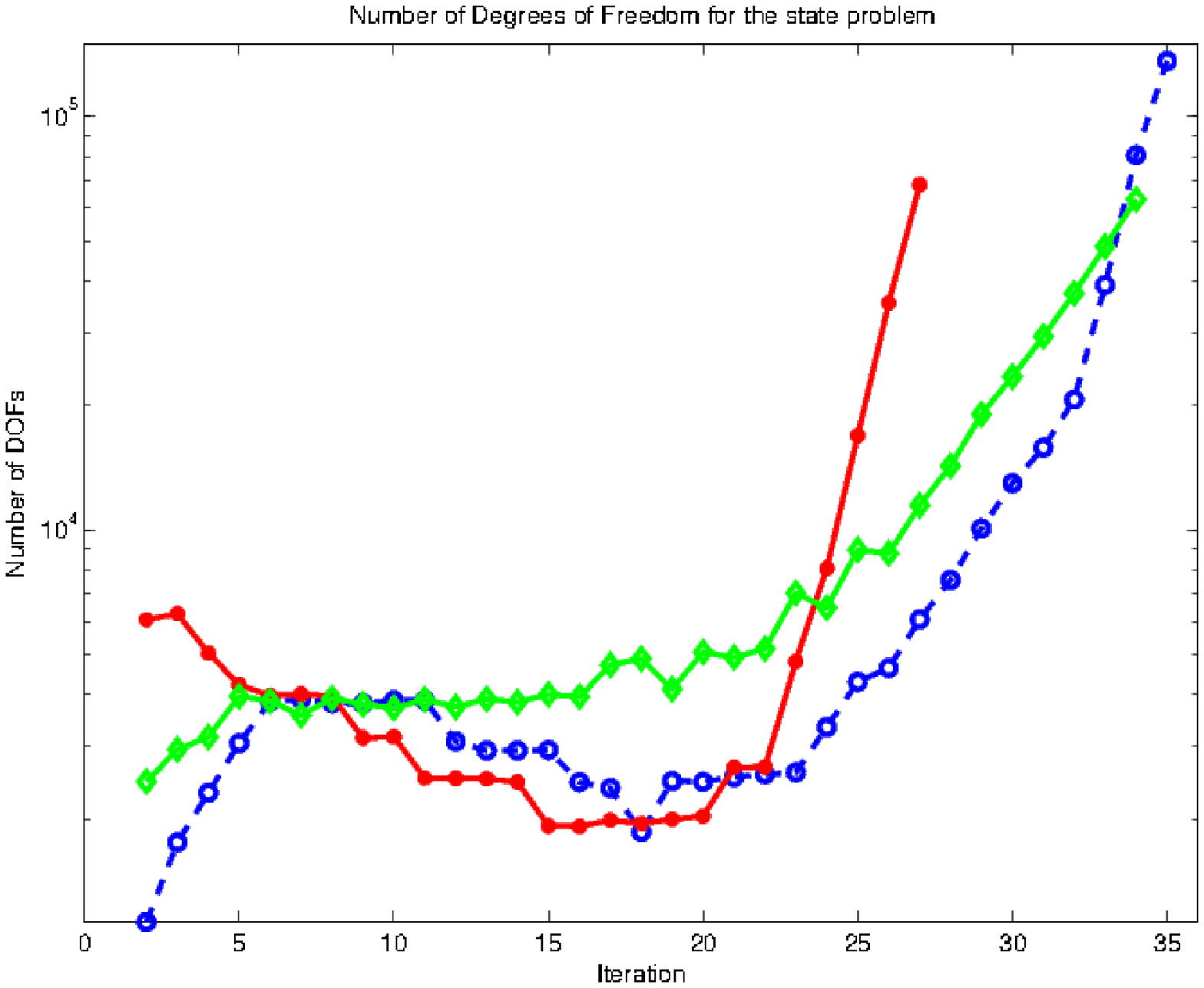}
    \label{fig:dofsEllInit}
    }
    \caption{Certified Descent Algorithm using multiple measurements for 
    initial guess (I) blue circles; (II) red stars; (III) green diamonds.
    (A) Evolution of the objective functional. (B) Number of Degrees of Freedom.}
    \label{fig:multiEllObjDof}    
\end{figure}

It is well-known in the literature that gradient-based strategies are local, that is 
they are able to identify only local minima. Thus, a key aspect 
for the success of the optimization procedure is represented by the initialization 
of the unknown variable. For shape optimization problems, this reduces to the 
choice of the initial shape. In figure \ref{fig:ellipseInterfaceInitializations}, 
we present three different initializations for the inclusion $\Omega$ and the 
corresponding final interfaces reconstructed by the Certified Descent Algorithm.
The three test cases confirm the ability of the method to correctly retrieve the 
portion of the inclusion close to $\partial\mathcal{D}$ whereas the regions in the 
center of the domain suffer from a degraded reconstruction.
Concerning the objective functional, the final values obtained using the 
proposed initializations are comparable (Fig. \ref{fig:objEllInit}), as the rapidly exploding number of Degrees 
of Freedom (Fig. \ref{fig:dofsEllInit}) which testifies again the ill-posedness of the inverse problem.

\begin{rmrk}
In the literature concerning inverse problems, a key issue when discussing a new 
method is represented by its robustness to noise and data perturbations. 
It is straightforward to observe that the construction of the sets of boundary data 
$(g,U_D)$ introduced at the beginning of this section is responsible for an additional 
contribution to the error. 
As previously observed, in order to retrieve reliable information to optimize the objective functional
the Certified Descent Algorithm requires an extremely high precision after few tens of iterations. 
In real-world applications, this results to be completely unfeasible since the additional information 
arising from the high precision would be lost due to the noisy nature of the boundary data.
Hence, for the purpose of this work we neglect the contribution of the uncertainty due to the 
measurements. 
An interesting extension of the current CDA may focus on the role of the error on the boundary measurements: 
in particular, an additional criterion may be integrated into the method in order to 
stop the algorithm if the error on the solutions of the state problems is smaller than the error 
on the data, that is if the effect of data fluctuations becomes predominant in the certification procedure.
\end{rmrk}


\subsection{The case of multiple boundary measurements}

\begin{figure}[p]
    \subfloat[$194$ elements.]
    {
    \includegraphics[width=0.3 \columnwidth]{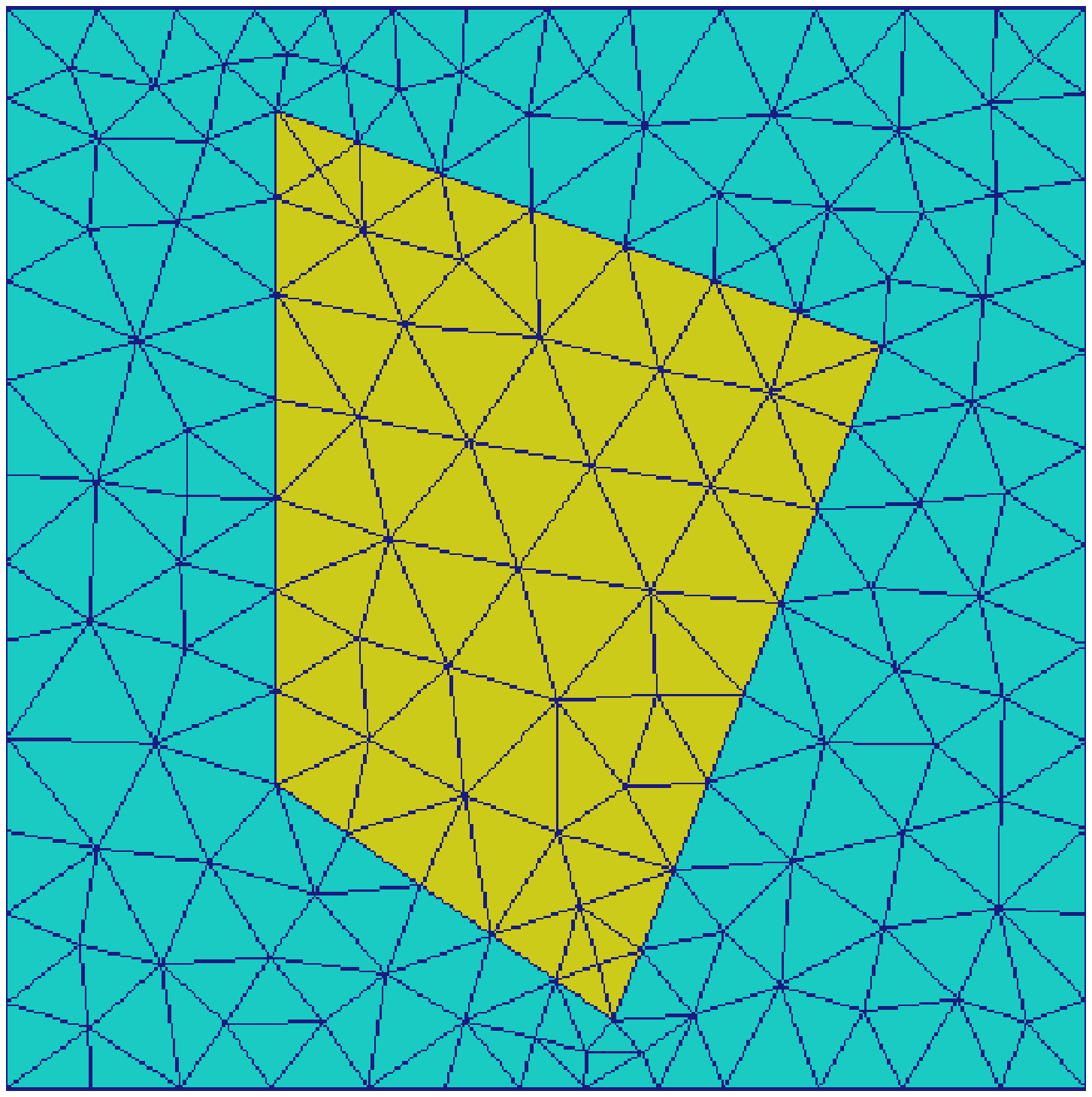}
    \label{fig:meshSquare1i}
    }
    \hfil
    \subfloat[$132040$ elements.]
    {
    \includegraphics[width=0.3 \columnwidth]{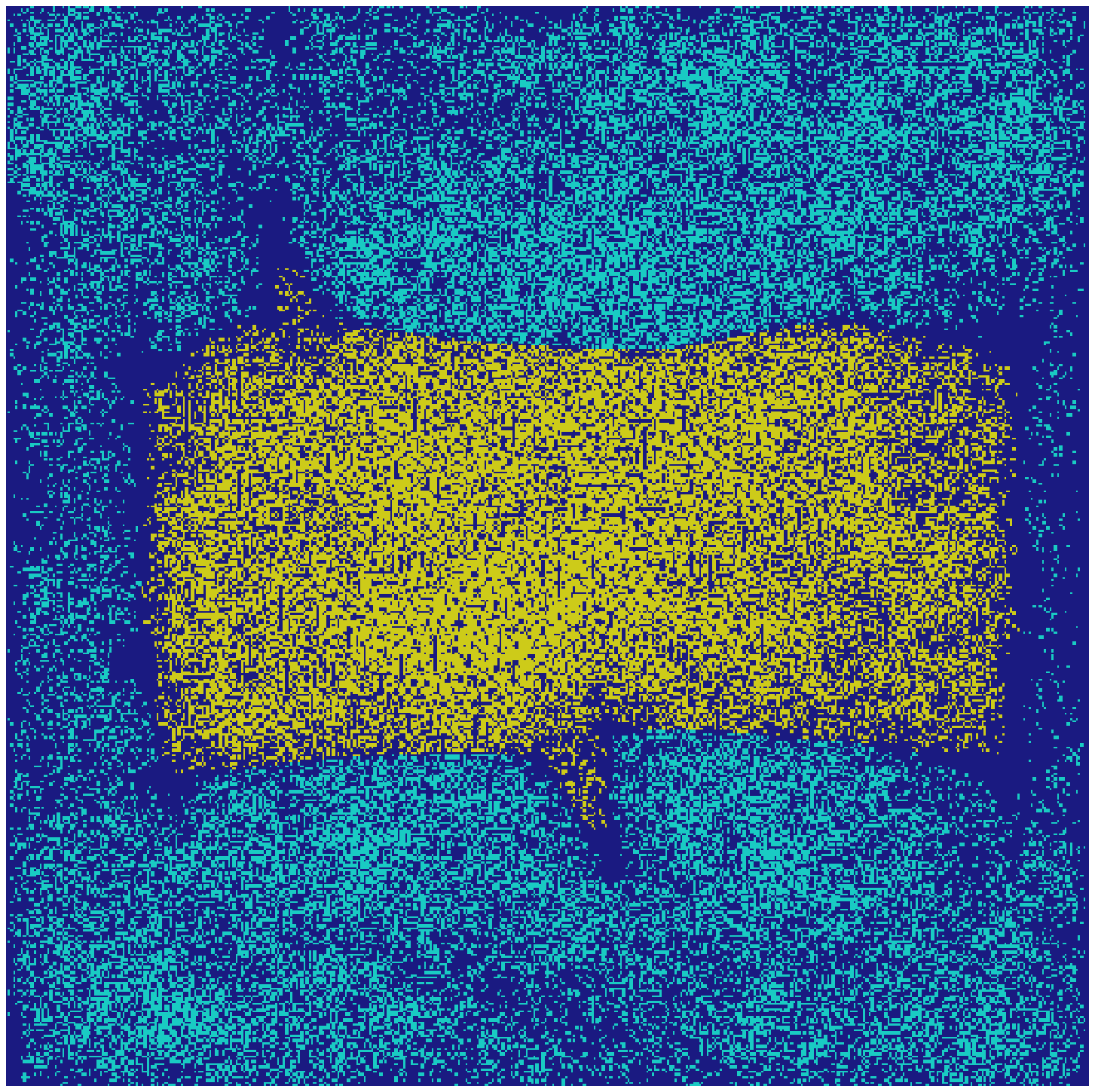}
    \label{fig:meshSquare1f}
    }
    \hfil
    \subfloat[Reconstructed interface.]
    {
    \includegraphics[width=0.3 \columnwidth]{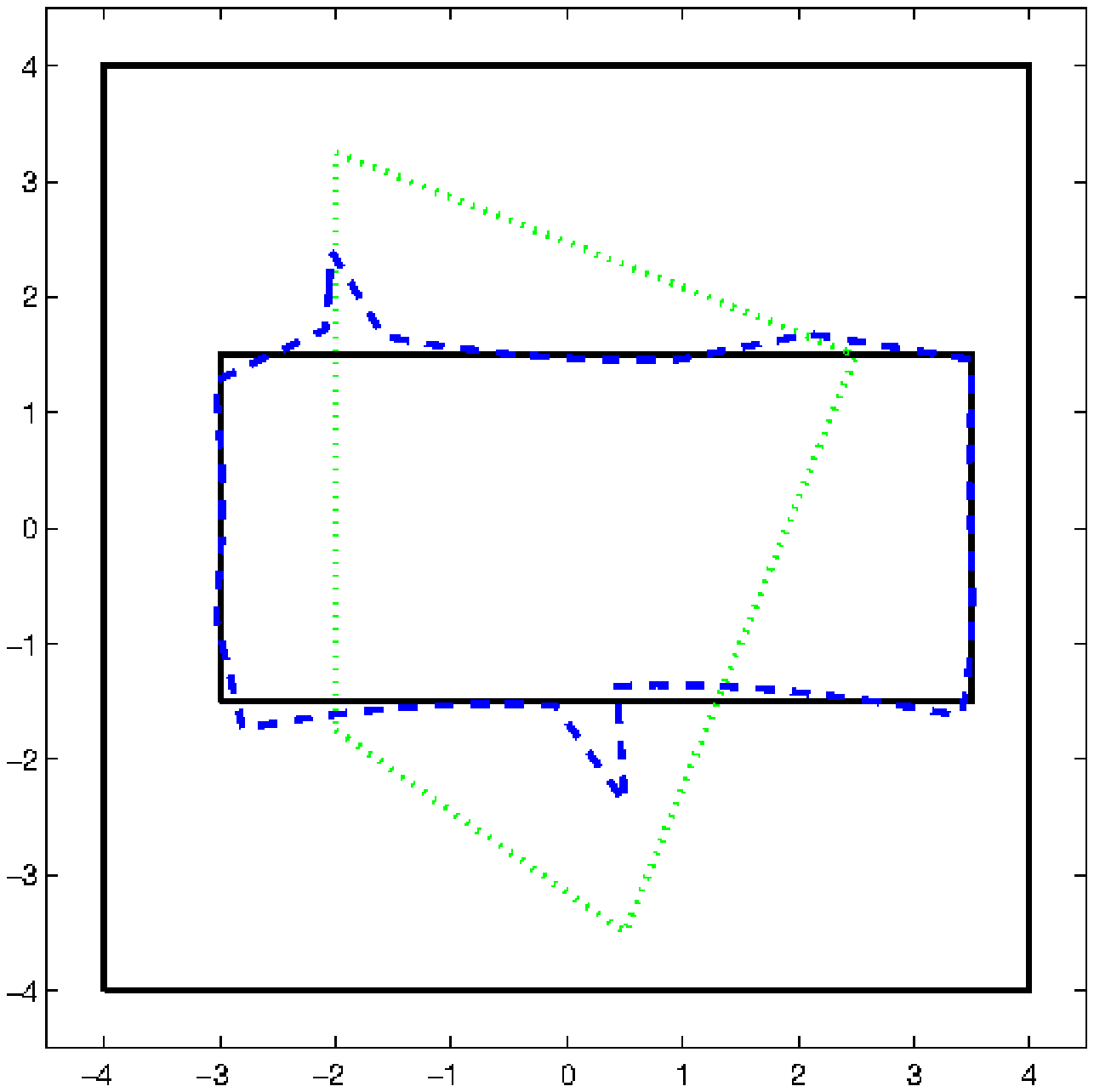}
    \label{fig:reconstructedSquare1}
    }

    \subfloat[$939$ elements.]
    {
    \includegraphics[width=0.3 \columnwidth]{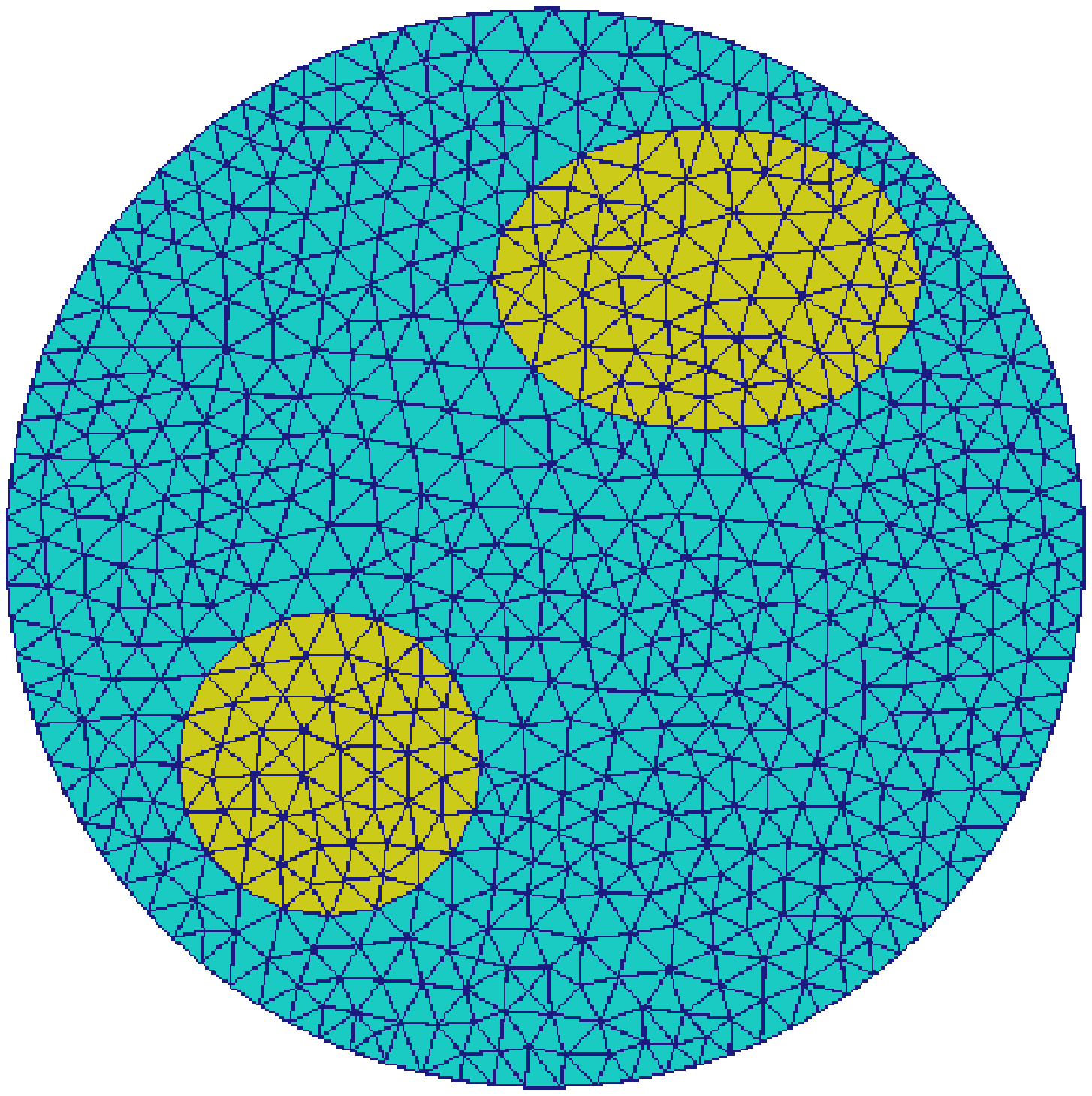}
    \label{fig:meshEllipse2i}
    }
    \hfil
    \subfloat[$428251$ elements.]
    {
    \includegraphics[width=0.3 \columnwidth]{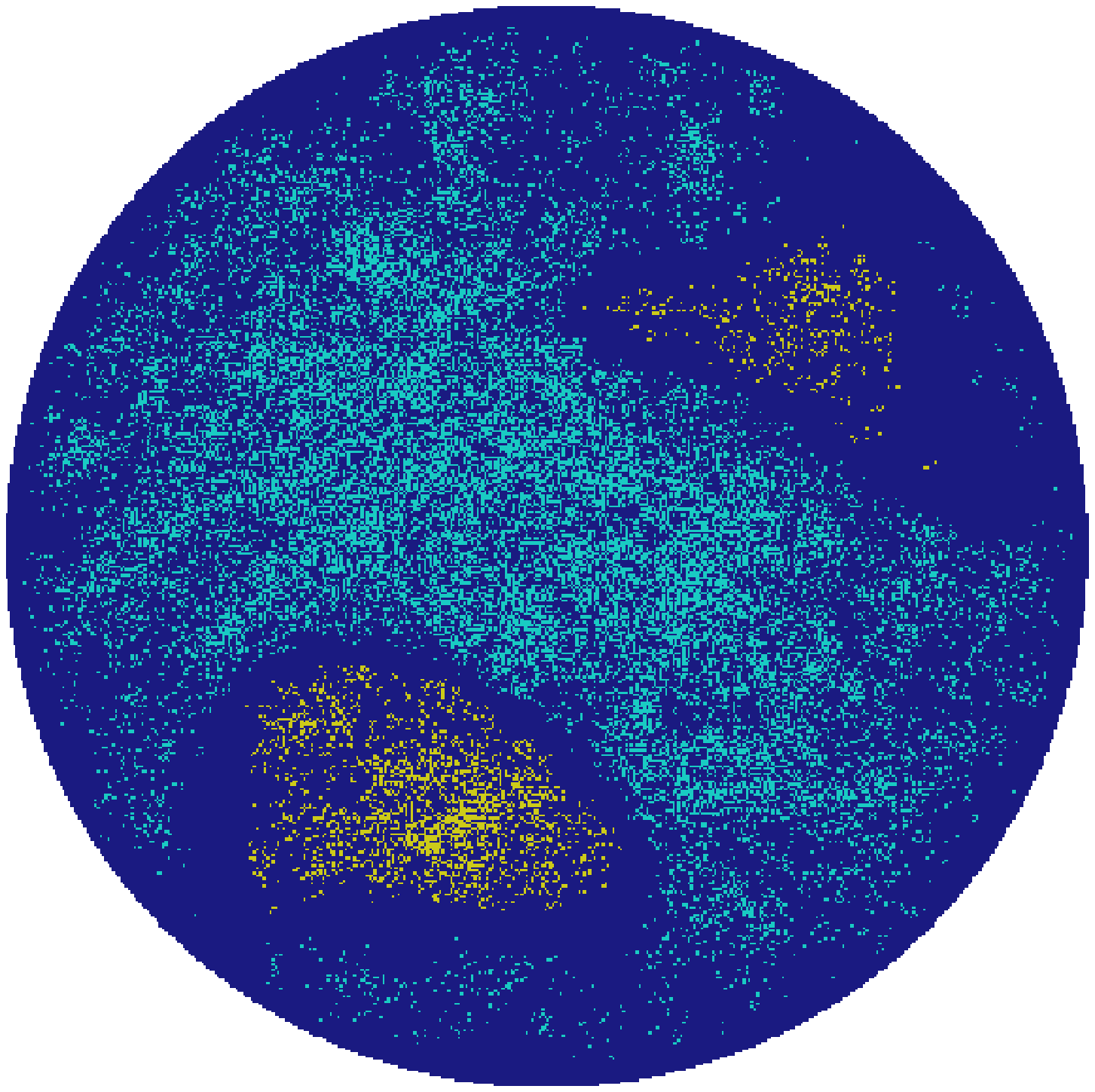}
    \label{fig:meshEllipse2f}
    }
    \hfil
    \subfloat[Reconstructed interface.]
    {
    \includegraphics[width=0.3 \columnwidth]{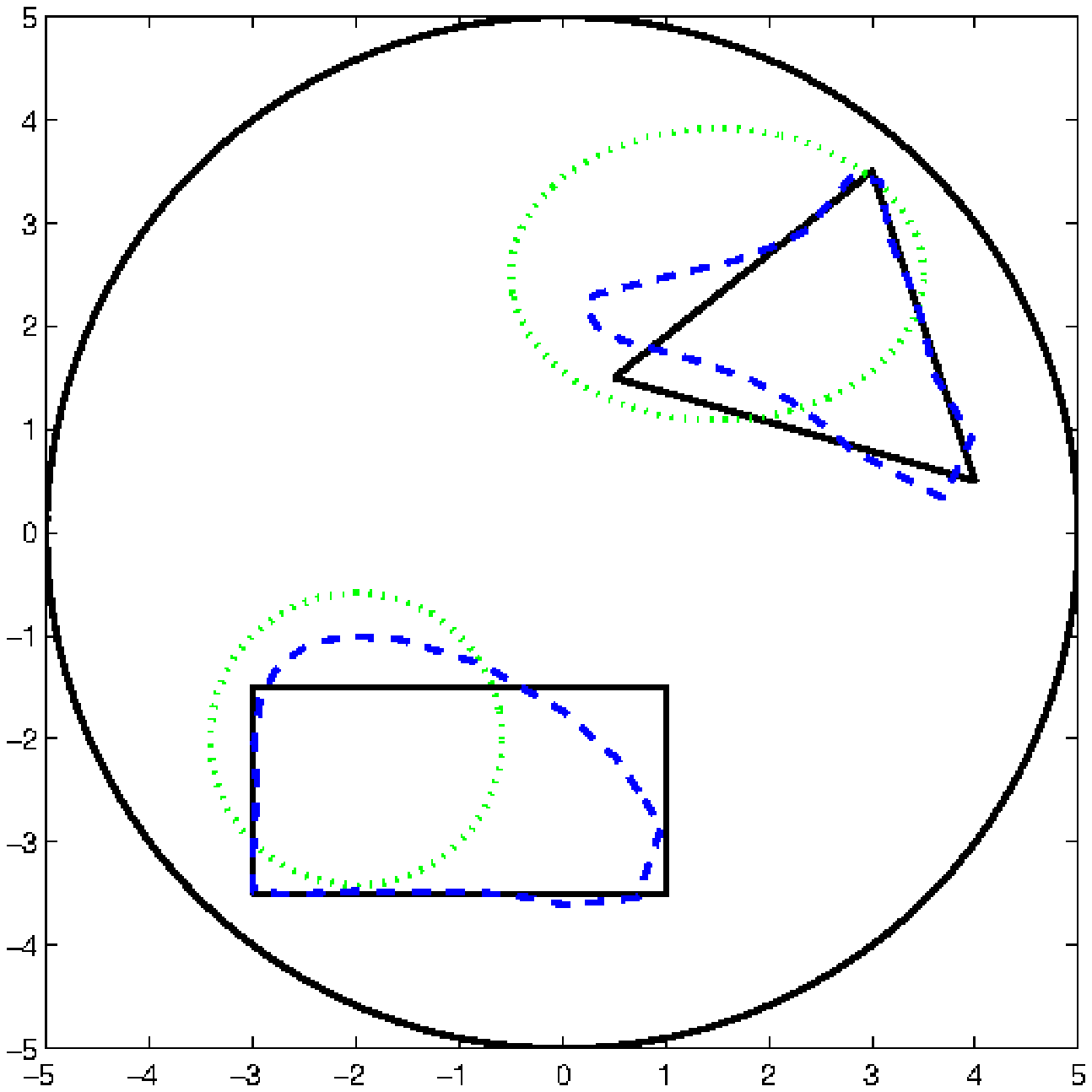}
    \label{fig:reconstructedEllipse2}
    }
    
    \subfloat[$1616$ elements.]
    {
    \includegraphics[width=0.3 \columnwidth]{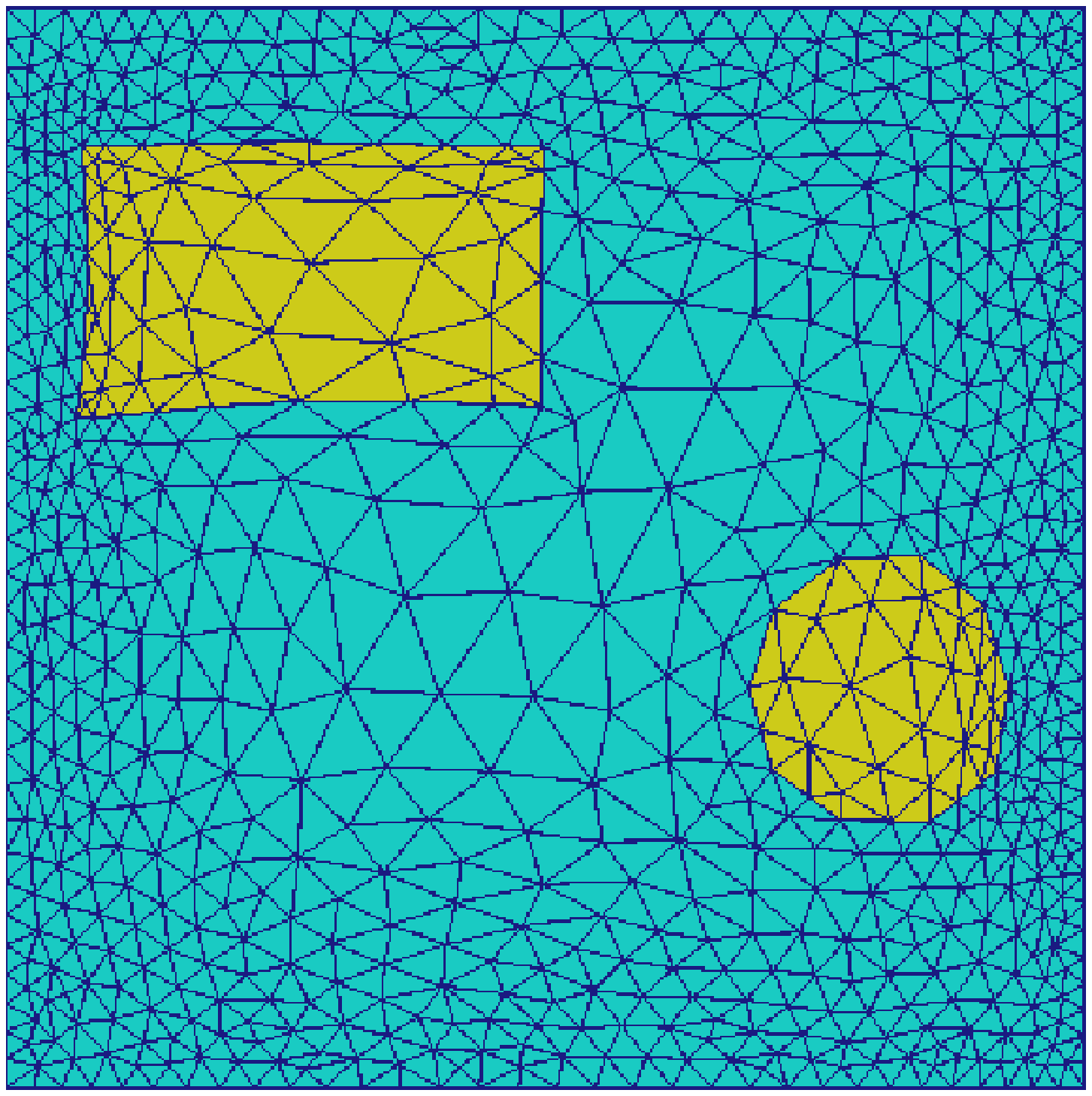}
    \label{fig:meshSquare2i}
    }
    \hfil
    \subfloat[$475744$ elements.]
    {
    \includegraphics[width=0.3 \columnwidth]{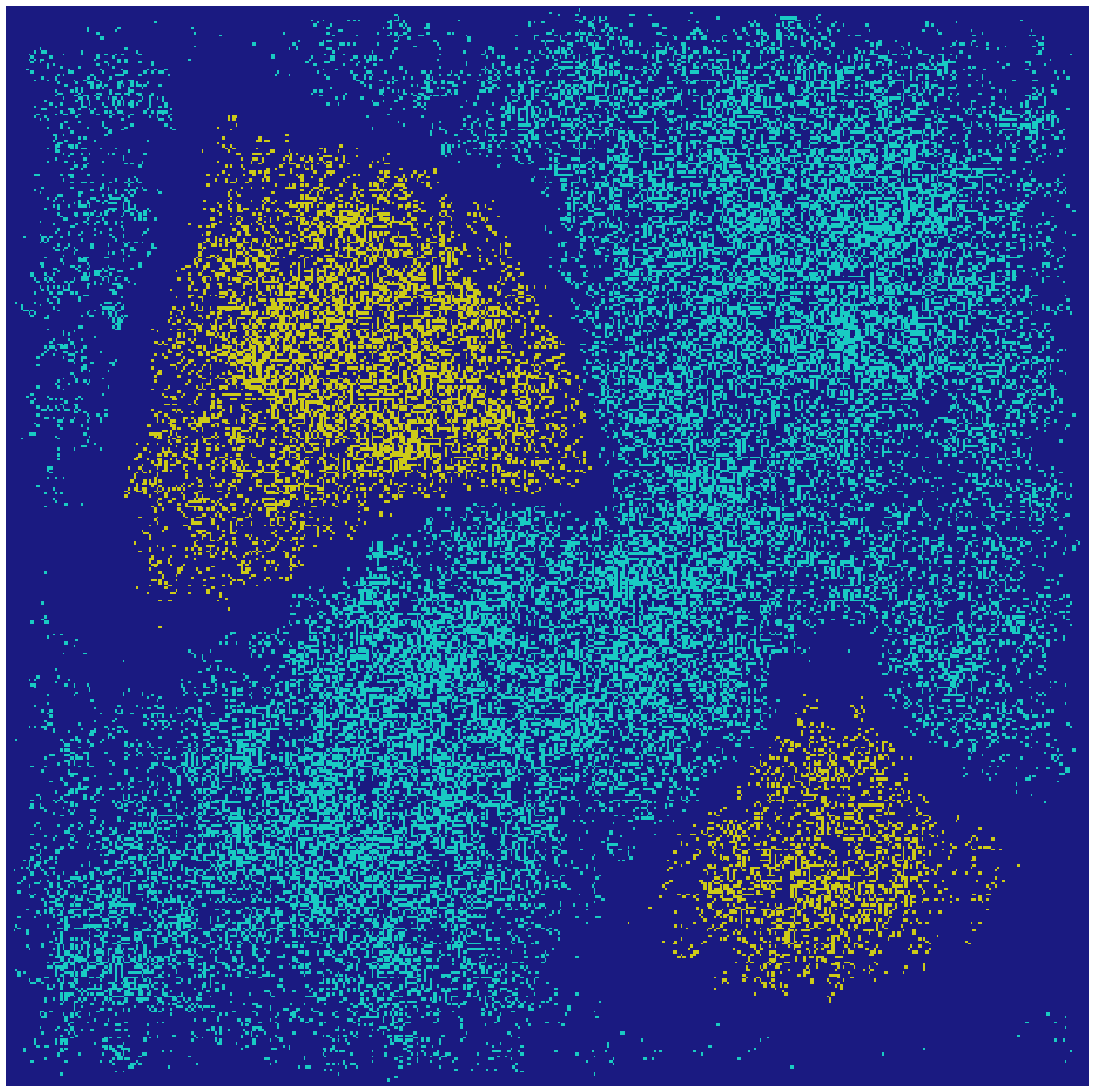}
    \label{fig:meshSquare2f}
    }
    \hfil
    \subfloat[Reconstructed interface.]
    {
    \includegraphics[width=0.3 \columnwidth]{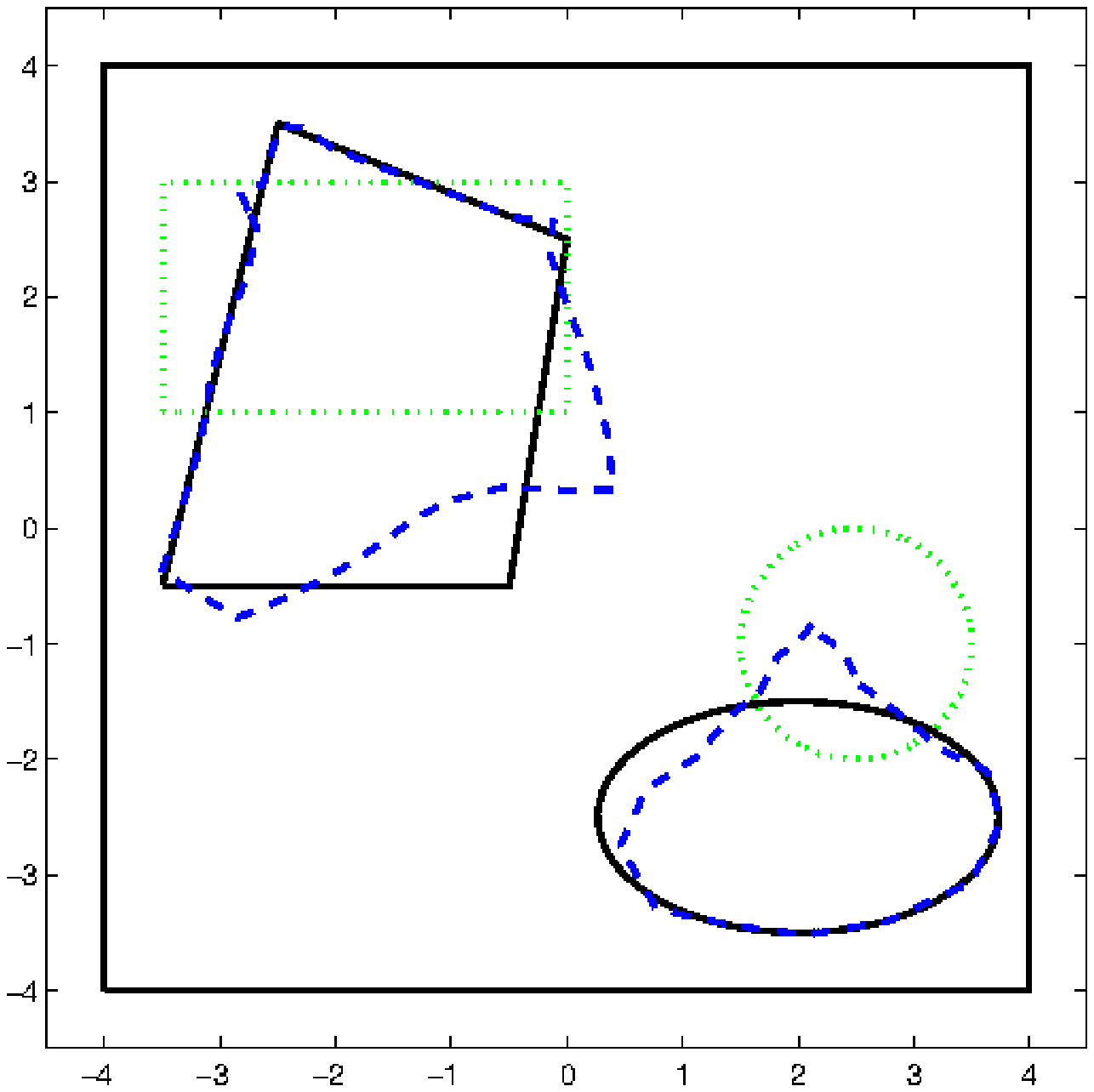}
    \label{fig:reconstructedSquare2}
    }
    \caption{Certified Descent Algorithm using multiple measurements: test cases 
(i)-(iii). Left: Initial mesh. Center: Final mesh. Right: Initial configuration 
(dotted green), target inclusion (solid black) and reconstructed interface 
(dashed blue). 
    }
    \label{fig:multiMeasures}    
\end{figure}

In \cite{2386935}, the authors prove that from an analytical point of view, one 
measurement is sufficient to uniquely reconstruct the inclusion within the 
Calder\'on's problem. Other analytical results on this topic are presented in 
\cite{doi:10.1080/01630560701381005}.
Nevertheless, from a numerical point of view, it is known that multiple 
measurements are required to have a correct approximation of the Electrical 
Impedance Tomography identification problem.
In this section, we present several tests of the previously described algorithm 
using multiple boundary measurements. In particular, we consider $D=10$ 
measurements such that $\forall j=0,\ldots,D-1$
$$
g_j(x,y) = (x+ a_j y)^{b_j} a_j^{c_j} \quad , \quad a_j=1+0.1j \quad ,  \quad 
b_j=\frac{j+1}{2} \quad , \quad c_j=j - 2 \left\lfloor \frac{j}{2} \right\rfloor
$$
and we use them to test the following cases:
\begin{enumerate}
\item[(i)] one inclusion in a square domain (Fig. 
\ref{fig:meshSquare1i}-\ref{fig:reconstructedSquare1});
\item[(ii)] two inclusions in a circular domain (Fig. 
\ref{fig:meshEllipse2i}-\ref{fig:reconstructedEllipse2});
\item[(iii)] two inclusions in a square domain (Fig. 
\ref{fig:meshSquare2i}-\ref{fig:reconstructedSquare2}).
\end{enumerate}
First, we present a simulation in which the body is the square $\mathcal{D} 
\coloneqq [-4,4]^2$ featuring a single polygonal inclusion (Fig. 
\ref{fig:reconstructedSquare1}).
Then, we propose two cases with multiple inclusions (Fig. 
\ref{fig:reconstructedEllipse2} and \ref{fig:reconstructedSquare2}): in both 
simulations, we assume that the number of inclusions is known \emph{a priori} 
and equals $2$ and that the conductivity $k_\Omega$ has only two values, one 
inside the inclusions and one for the background.

As expected, the use of multiple measurements provides sufficient information to 
reconstruct the target interfaces near the boundary $\partial\mathcal{D}$ (Fig. 
\ref{fig:reconstructedSquare1}, \ref{fig:reconstructedEllipse2} and 
\ref{fig:reconstructedSquare2}). 
On the contrary, the identification of the inner interfaces appears more 
difficult and with a less precise outcome.
This phenomenon is due to the severe ill-posedness of the problem and the 
diffusive nature of the state equations is responsible for the loss of the 
information in the center of $\mathcal{D}$, even when using several 
measurements. 
In particular, we remark that the final interface in figure 
\ref{fig:reconstructedSquare1} still presents two kinks from the initial 
configuration: this is mainly due to the aforementioned phenomenon and may be 
bypassed by choosing a regularizing scalar product. However, this approach would 
lead to a global smoothing of the reconstructed interface, including a loss of 
information about the potential sharp physical corners of the polygonal 
inclusion. 

\begin{figure}[htb]
\centering
    \subfloat[Objective functional.]
    {
    \includegraphics[width=0.48\columnwidth]{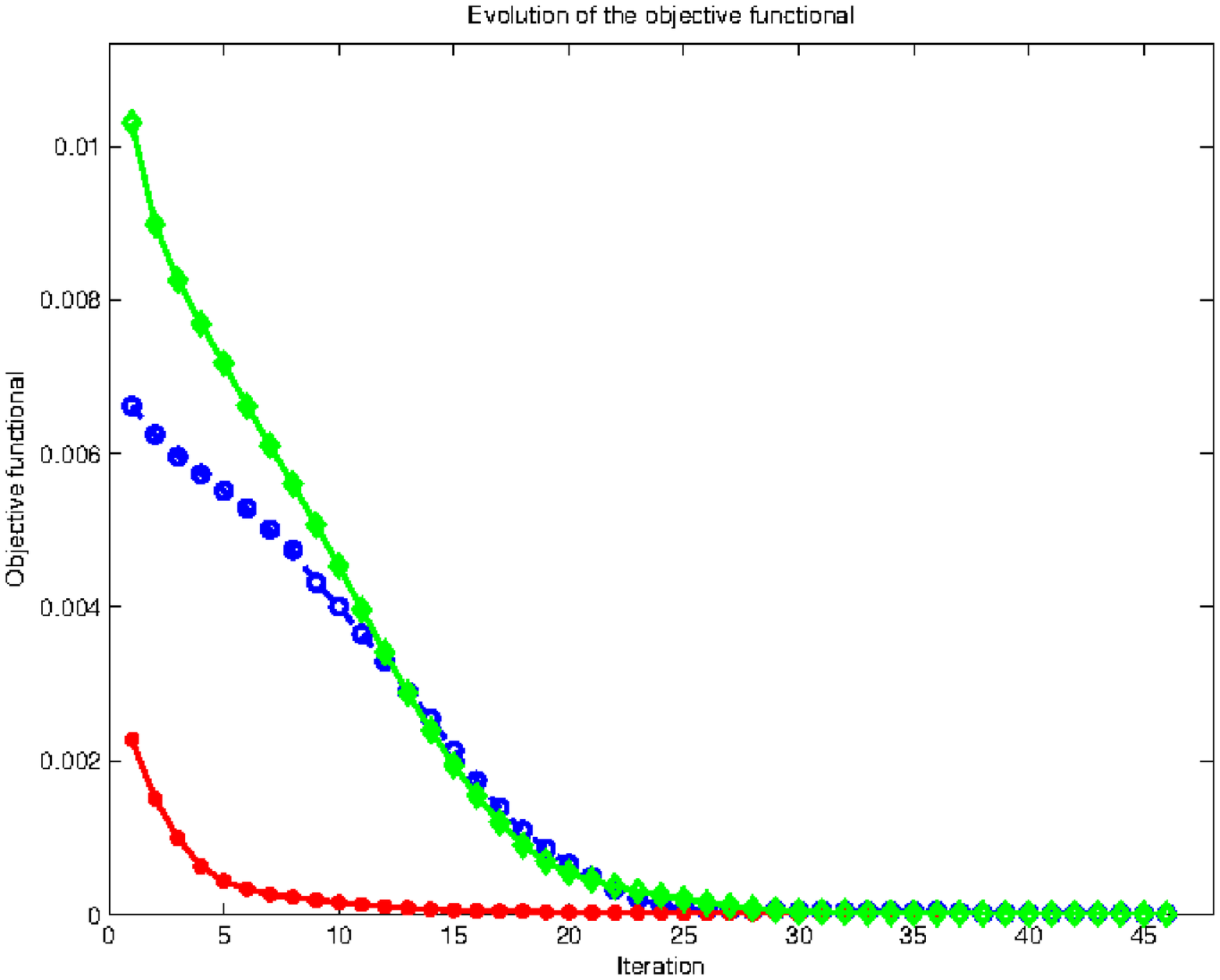}
    \label{fig:objMulti}
    }
    \hfil
    \subfloat[Number of Degrees of Freedom.]
    {
    \includegraphics[width=0.48\columnwidth]{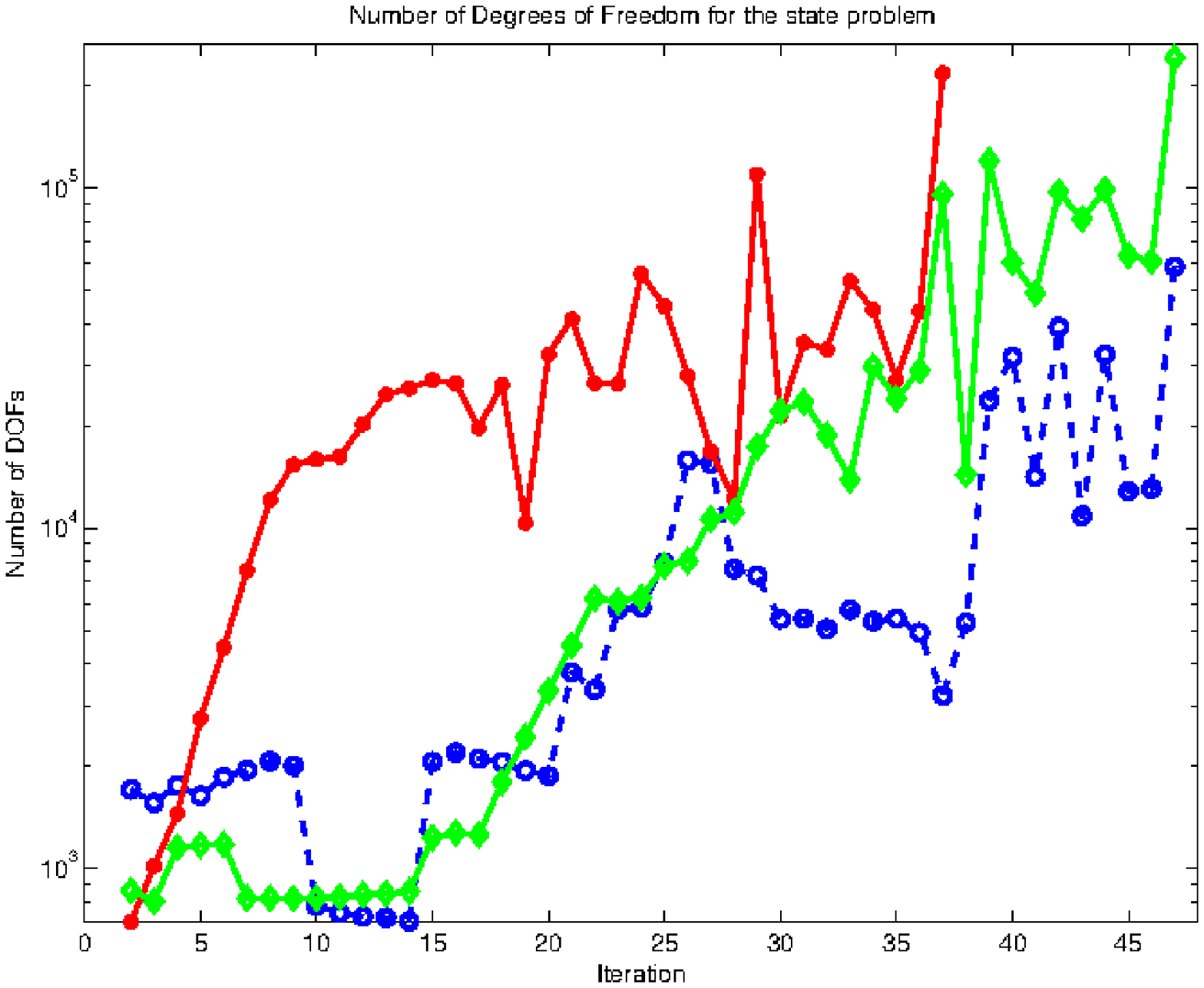}
    \label{fig:dofsMulti}
    }
    \caption{Certified Descent Algorithm using multiple measurements: test case 
(i) blue circles; (ii) red stars; (iii) green diamonds. (A) 
Evolution of the objective functional. (B) Number of Degrees of Freedom.}
    \label{fig:multiObjDof}    
\end{figure}

As previously remarked, figure \ref{fig:objMulti} confirms the monotonically 
decreasing behavior of the objective functional with respect to the iterations 
of the algorithm. Moreover, the quantitative information associated with the 
estimator of the error in the shape gradient allows to derive a reliable 
stopping criterion for the optimization procedure which results to be 
fully-automatic. \\
Eventually, coarse meshes are proved to be reliable for the computation during 
the initial iterations when the guessed position and shape of the inclusion is 
very unlikely to be precise. 
Within this context, even few Degrees of Freedom provide enough information to 
identify a genuine descent direction for the objective functional which we 
later certify using the discussed goal-oriented procedure. Thus, the same meshes 
may be used for several iterations increasing the number of Degrees of Freedom 
only when the descent direction is no more validated (Fig. 
\ref{fig:dofsMulti}). 
\\
Both the inability of the method to reconstruct the interface far from the external boundary 
and the rapidly increasing number of Degrees of Freedom required to certify the descent 
direction clearly testifies the limitations of classical gradient-based approaches when 
dealing with the problem of Electrical Impedance Tomography.
The ill-posedness of the problem is confirmed by the fact that after few tens of 
iterations, we are unable to identify a genuine descent direction at a small 
computational cost since the number of Degrees of Freedom required by the certification 
procedure rapidly reaches $10^5$.
Nevertheless, the main novelty of the Certified Descent Algorithm - that is its certification 
procedure - provides us an heuristic criterion to stop the optimization routine when the 
number of Degrees of Freedom tends to explode, being the 
improvement of the solution negligible with respect to the huge precision the computation 
would require.

\section{Conclusion}
\label{ref:conclusion}

In this work, we coupled classical shape optimization techniques with a 
goal-oriented error estimator for the shape gradient. 
A guaranteed bound for the error in the shape gradient has been derived by means 
of a certified \emph{a posteriori} estimator. \\
On the one hand, introducing an \emph{a posteriori} estimator for the error in 
the shape gradient provides quantitative information to define a reliable 
stopping criterion for the overall optimization procedure. 
Coupling this approach with the 2-mesh shape optimization strategy introduced in 
\cite{smo-AP} results in the novel Certified Descent Algorithm. The CDA is a 
fully-automatic procedure for certified shape optimization: a validation of the 
method is presented by means of several test cases for the well-known inverse 
identification problem of Electrical Impedance Tomography. 
Another important feature of this method is the ability of identifying a 
certified descent direction at each iteration thus leading to a monotonically 
decreasing evolution of the objective functional.

Even though the CDA is able to make coarse meshes reliable to identify a genuine 
descent direction for the objective functional during the initial iterations, 
the overall computational cost tends to remain high. 
As a matter of fact, the major drawback of the described procedure is the 
necessity of solving the dual flux problems to derive a fully-computable upper 
bound of the error in the shape gradient. 
Hence, the Certified Descent Algorithm may result in higher computing times than 
the Boundary Variation Algorithm applied on fine meshes.

Ongoing research focuses on improving the \emph{a posteriori} estimates for the 
discretization error in the shape gradient. 
Promising results are expected by the development of error estimators that only 
involve the computation of local quantities. 
Within this framework, accounting for anisotropic mesh adaptation 
\cite{FLD:FLD2688} may lead to discretizations with a lower number of Degrees of 
Freedom and a better approximation of the physical problem. \\
Future investigations will focus on the application of the Certified Descent 
Algorithm to other challenging problems, such as the shape optimization of 
elastic structures.

\bibliographystyle{abbrv}
\bibliography{./COCV_art2015-92_bibliography}

\end{document}